\documentclass{amsart}

\usepackage{amsthm,amssymb,verbatim}
\usepackage{graphicx}
\usepackage{enumerate}

 \newcommand{\Aa}{\mathcal A} %space of annuli
 \newcommand{\Cc}{\mathcal C}
 
 \newcommand{\Nn}{\mathcal{N}}

 \newcommand{\FF}{\mathcal{F}}
 \newcommand{\RR}{\mathbf{R}}  % reals
 \newcommand{\ZZ}{\mathbf{Z}}  % integers
 \newcommand{\BB}{\mathbf{B}}  %ball
 \newcommand{\CC}{\mathbf{C}}  %complex numbers
 \newcommand{\dist}{\operatorname{dist}}

 \newcommand{\eps}{\epsilon}
 \newcommand{\Tan}{\operatorname{Tan}}
 \newcommand{\Length}{\operatorname{Length}}

\input epsf
\def\begfig {
\begin{figure}
\small }
\def\endfig {
\normalsize
\end{figure}
}

\swapnumbers

    \newtheorem{theorem}    {Theorem}       [section]
    \newtheorem{lemma}      [theorem]       {Lemma}
    \newtheorem{corollary}  [theorem]     {Corollary}
    \newtheorem{proposition}       [theorem]       {Proposition}
    
    \newtheorem*{claim}{Claim}
    \newtheorem*{theorem*}{Theorem}

    \theoremstyle{definition}
    \newtheorem{definition}  [theorem] {Definition}

    \theoremstyle{definition}
    \newtheorem{remark}   [theorem]       {Remark}

\begin{document}

\renewcommand{\thesubsection}{\thetheorem}
   % this is so subsections and theorems, etc will be
   % numbered together.

\title[Genus-one helicoids]{Genus-one helicoids \\from a variational point of view}
\author{David Hoffman}
\thanks{The research of the first author was supported by the NSF under
   grant~DMS-0139410}
\address{Department of Mathematics\\ Stanford University\\ Stanford, CA 94305}
\email{hoffman@math.stanford.edu}
\author{Brian White}
\thanks{The research of the second author was supported by the NSF
  under grants~DMS-0406209 and~DMS-0707126}
\email{white@math.stanford.edu}
\date{November 20, 2006. Revised January 24, 2008.}

\maketitle

\section{Introduction}\label{section:intro}

In this paper, we prove by variational means the existence of a
complete, properly embedded, genus-one minimal surface in $\RR^3$
that is asymptotic to a helicoid at infinity.   We also prove
existence of surfaces that are asymptotic to a helicoid away from
the helicoid's axis, but that have infinitely many handles
arranged periodically along the axis. These theorems were
originally proved by very different methods in~\cite{wwh1}. We
also prove some new properties of such helicoid-like surfaces.

To state the theorems precisely, we need some terminology. For
$s\in \RR$ we let $\sigma_s$  denote the screw motion of $\RR^3$
defined by
\begin{equation*}
   \sigma_s(\cos\theta, \sin\theta, z)
     =
   (\cos(\theta+s), \sin(\theta+s), z+s ).
\end{equation*}
We let $H$ be the standard helicoid that contains the $x$-axis $X$
and the $z$-axis $Z$, and that is invariant under the screw
motions $\sigma_s$.

A {\em nonperiodic genus-one helicoid} is a complete, properly
immersed minimal surface in $\RR^3$ that is conformally a
once-punctured torus (i.e., a torus with one point removed) and
that is asymptotic to $H$ at infinity. If it also contains $X$ and
$Z$, we say that it is {\em symmetric}: by Schwarz reflection, it
is invariant under the $180^\circ$ rotations $\rho_X$ and $\rho_Z$
about $X$ and $Z$, and hence it is also invariant under their
composition $\rho_Y$.  (Note that the symmetry group of any
positive-genus surface asymptotic to $H$ must be a finite subgroup
of the symmetry group $\operatorname{sym}(H)$ of $H$.  It is not
hard to prove that any finite subgroup of $\operatorname{sym}(H)$
is conjugate to a subgroup of
  $\{ I, \rho_X, \rho_Y, \rho_Z\}$,  and thus that a symmetric genus-one helicoid is as symmetric as a genus-one helicoid can possibly be.)

\begin{theorem}\label{thm:first}
There exists an embedded, symmetric, nonperiodic genus-one
helicoid.
\end{theorem}

A {\em periodic genus-one helicoid} is a complete, properly
immersed minimal surface $N$ in $\RR^3$ such that for some $h>0$:
\begin{enumerate}
 \item $N$ is $\sigma_{2h}$-invariant,
 \item $N/\sigma_{2h}$ is conformally a twice-punctured torus, and
 \item $N/\sigma_{2h}$ is asymptotic to $H/\sigma_{2h}$ at
 infinity.
\end{enumerate}
If, in addition, $N$ contains $X$ and $Z$ and if the fundamental
domain
$$
    N \cap \{(x,y,z): |z|<h\}
$$
is bounded by the two lines $\sigma_{\pm h}(X)$, then $N$ is
called a {\em symmetric} periodic genus-one helicoid. (Note that
the quotient $N/\sigma_{2h}$ of a symmetric periodic genus-one
helicoid $H$ contains two horizontal lines, corresponding to the
two parities of
  $n$
in the lines $\sigma_{nh}(X)$.)

If $N$ is any genus-one helicoid, let $h(N)$ be the smallest $h>0$
such that $N$ is $\sigma_{2h}$-invariant.  (If $N$ is nonperiodic,
we let $h(N)=\infty$.)

\begin{theorem}\label{thm:second}
For every $h\in (\pi/2,\infty]$, there is an embedded, symmetric,
genus-one helicoid $N$ with $h(N)=h$.
\end{theorem}

Note that Theorem~\ref{thm:first} is the special case $h=\infty$
of Theorem~\ref{thm:second}.

The condition $h>\pi/2$ is sharp: in~\cite{hoffwhite2}, we prove
that there are no examples with $h(N)\le \pi/2$, even if we allow
somewhat less symmetry.  (For $h(N)<\pi/2$ this was observed by
Meeks.)

Uniqueness is known for $h=\pi$ \cite{wwh1, fm05}, but not in
general. Thus there may be other embedded, symmetric genus-one
helicoids that do not arise from the construction in this paper.
The following theorem holds for {\em all} symmetric genus-one
helicoids:

\begin{theorem}\label{thm:third}
Let $\eta>\pi/2$. Let $N$ be an embedded, symmetric genus-one
helicoid with $h(N)=h\in [\eta, \infty]$, and let $M$ be the
fundamental domain
$$
  M = N \cap \{ |z| < h\}.
$$
\begin{enumerate}
 \item The intersection $N\cap H$ consists of $Z$ together with
 the horizontal lines
 $$
    \{(x,y,z)\in H: z = n h \}, \qquad (n\in \ZZ).
 $$
 Furthermore, $N\setminus H$ consists of congruent, simply
 connected components.  The fundamental domain $M$ consists of two of
 those components, one on either side of $H$.
 \item For each vertical plane $V$, there are at most four points
 $p\in M\setminus Z$ for which $\Tan_pN$ is parallel to $V$.
 Such points must lie in the cylinder
 $$
    \{(x,y,z): x^2+y^2\le R^2,\, |z|< 2\pi \}
 $$
 where $R$ depends only on $\eta$.
 \item The space of all such $N$ (for a given $\eta>\pi/2$) is
 compact with respect to smooth convergence on bounded subsets of
 $\RR^3$.
\end{enumerate}
\end{theorem}

If $h(N)=\infty$, then the fundamental domain $M$ is all of $N$,
and the horizontal lines in assertion~(1) of the theorem consist
just of $X$.

Concerning assertion~(1), note that by our definition of
``symmetric'', $N\cap H$ must contain the indicated lines.  It is
perhaps surprising that $N\cap H$ contains no other points.  For
the surfaces we produce, this property follows immediately from
the construction.  The proof that~(1) holds in general will be
given elsewhere~\cite{hoffwhite2}. Proofs of~(2) and~(3) are given
in Section~\ref{section: compactness}.

A proof of assertion~(3), the compactness result, by different
methods is implicit in~\cite{wwh1}.

  \begfig
 \hspace{0.3in}
  \vspace{.2in} \centerline {
  \includegraphics[width=2.5in]{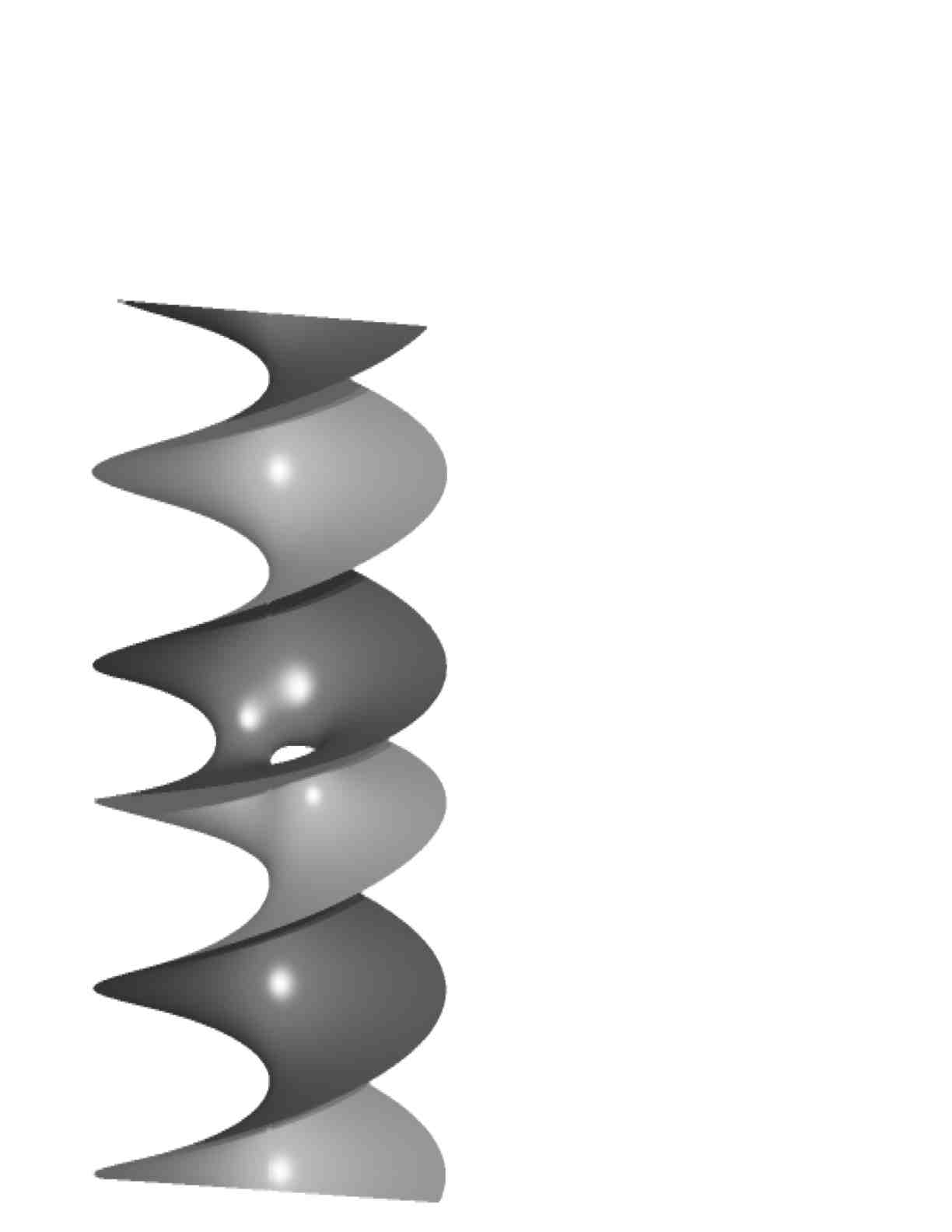}
  \includegraphics[width=2.05in]{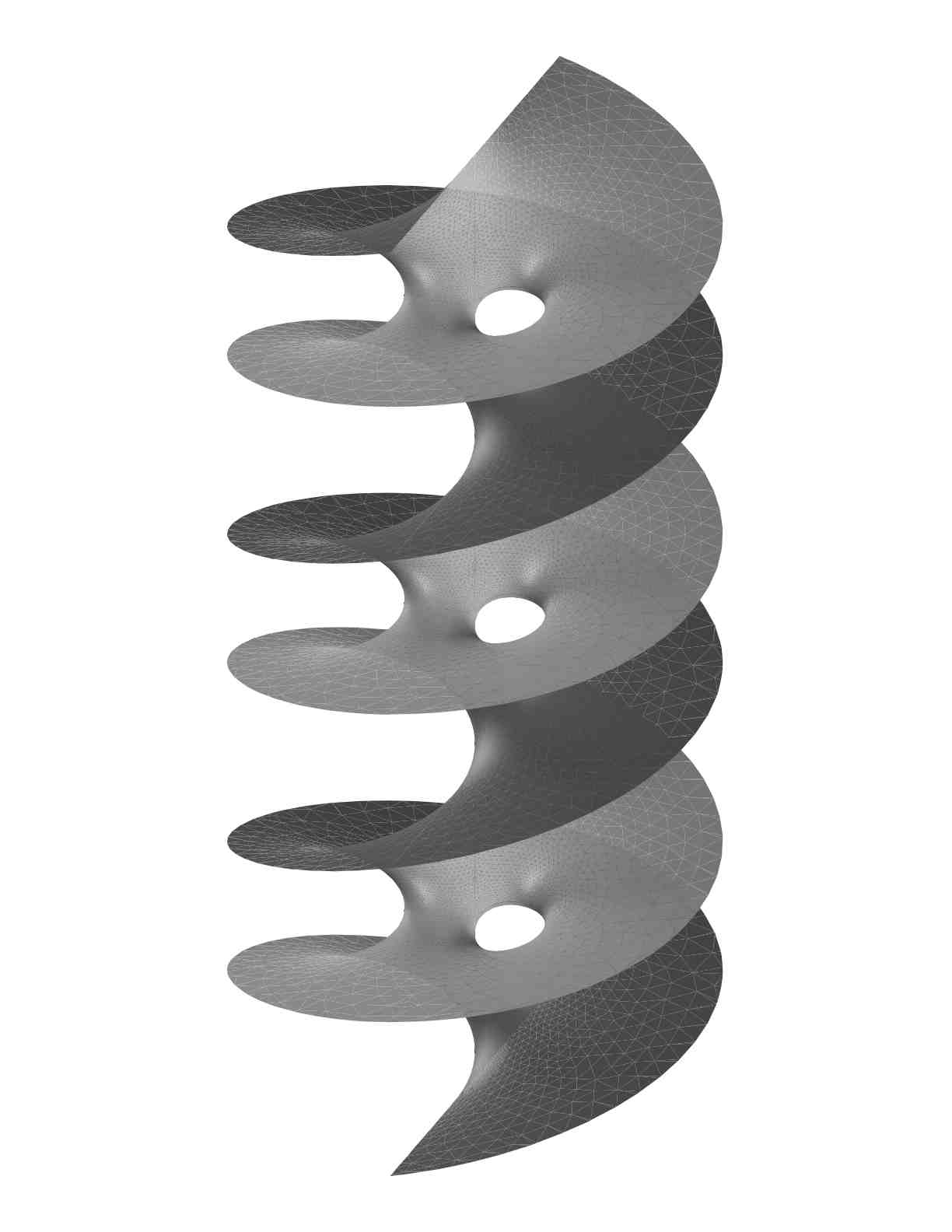}
   }
 \hfill
 \begin{center}
 \vspace{0.3cm}
 \parbox{4.1in}{
  \caption{\label{figure:genus-one helicoids} Nonperiodic (left) and periodic (right)
  genus-one helicoids.   }
 }
 \end{center}
 \endfig

Hoffman, Karcher, and Wei discovered symmetric genus-one
helicoids. (Whether asymmetric examples exist is not known.) In
\cite{howe3}, they proved existence but not embeddedness of an
example with $h=\infty$.  In \cite{hkw8}, they proved existence of
an example with $h=\pi$, and they proved that every $h=\pi$
example must be embedded.  By numerical computations, they also
discovered embedded examples for a range of values of $h(N)$; that
work is described in \cite{howe5}.

Hoffman, Weber, and Wolf \cite{wwh1} proved that there exists a
continuous family of embedded genus-one helicoids taking all
values of $h$ in $(\pi/2,\infty]$.  They also proved that there is
only one example with $h=\pi$.  This uniqueness was proved
independently in a different way by Ferrer and Martin~\cite{fm05}.

All of those investigations relied on the Weierstrass
representation, which requires solution of period problems and a
separate proof of embeddedness.  We present here a method of
realizing examples as limits of compact, embedded minimal
surfaces. Period problems do not arise, and the method gives
existence and embeddedness at the same time.  It also provides
additional information about the geometric behavior of the
surfaces.  Our investigation is similar in spirit to \cite{hm9},
where a variational construction of the generalized, higher-genus,
Costa surfaces with three ends is given, and to the part of
\cite{chm2} in which translation-invariant, Callahan-Hoffman-Meeks
surfaces of odd genus are produced.

In our view, the construction presented here gives a good answer
to the question of {\em why} genus-one helicoids should exist.

Using deep results of Colding and Minicozzi
\cite{ColdingMinicozzi1, ColdingMinicozzi2,
  ColdingMinicozzi3,
  ColdingMinicozzi4},
Meeks and Rosenberg proved \cite{MR01} that the helicoid is the
unique properly embedded, simply connected, nonplanar minimal
surface---there is no symmetry assumption in their result. Whether
or not there is only one nonperiodic, embedded, symmetric
genus-one helicoid---and more generally whether there exist any
nonsymmetric examples---has yet to be resolved.

The authors would like to thank Mike Wolf and Matthias Weber for
helpful conversations.

\section{Outline of the construction}
 \label{section:outline}

\stepcounter{theorem}
\subsection{The basic idea}\label{helicoidguide}

Let $H_{R,h}$ be the portion of $H$ inside a solid cylinder
centered at the origin, with axis $Z$, radius $R$, and height
$2h$:
\begin{equation*}
   H_{R,h} = H\cap\{|z|< h, \, r<R\}.
\end{equation*}
The boundary $\partial H_{R,h}$ is a simple closed curve
consisting of two horizontal line segments
   $L_{\pm h}\subset\sigma_{\pm h}(X)$
and two helical curves connecting the endpoints of those line
segments.  Note that $H_{R,h}$, like $H$, is invariant under
$180^\circ$ rotations about the coordinate axes.

Our idea is simple: replace $H_{R,h}$ by an embedded, genus-one
minimal surface $M_{R,h}$ with the same boundary.  Since we want
our genus-one helicoids to be symmetric, we require that $M_{R,h}$
contain the segments $X\cap H_{R,h}$ and $Z\cap H_{R,h}$.  By
taking a limit as $R$ and $h$ tend to infinity, we hope to get a
nonperiodic genus-one helicoid. By taking a limit as $R\to\infty$
with $h$ fixed, we hope to get a fundamental domain for a
$\sigma_{2h}$-invariant periodic genus-one helicoid.

 \begfig
 \hspace{0.3in}
\vspace{.2in} \centerline {
\includegraphics[width=3in]{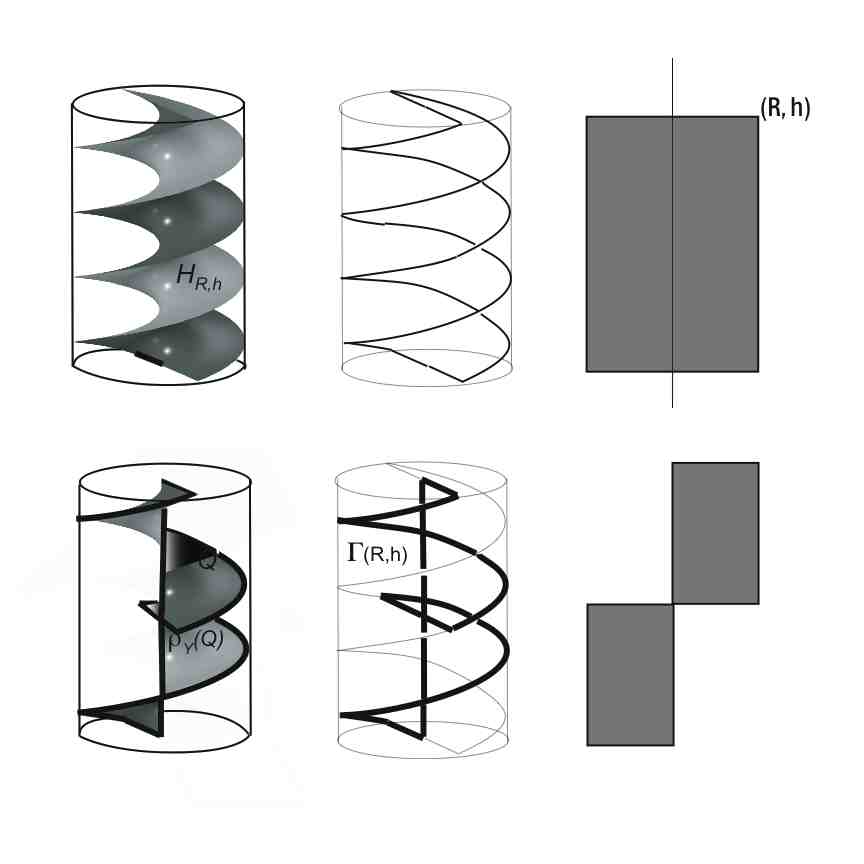}
} \vspace{-.5in}
 \hfill
 \begin{center}
 \vspace{0.3cm}
  \caption{} \label{figure:sixpics}\vspace{.3cm}
   \parbox{4.0in}{
 Top left: The surface $H_{R,h}= H\cap C$, where $C$
is the solid cylinder centered at the origin with axis $Z$, height
$2h$ and radius $R$. Top center: $\partial H_{R,h}= H\cap\partial
C$. Top right: The parameter domain of $H_{R,h}$, assuming the
parametrization
   $(u,v)\to (u\cos v, u\sin v ,v)$.
 \newline Bottom left: The union $Q_{R,h}$ of two of the four congruent
components of $H_{R,h}\setminus(X\cup Z)$.  Bottom right: The
corresponding region in the parameter domain.  Bottom center:
$\Gamma(R,h):=\partial Q_{R,h}$.
  }
 \end{center}
 \endfig

\stepcounter{theorem}
\subsection{Replacing $H_{R,h}$ by a genus-one surface}
  \label{subsection:replacing}
If $X$ and $Z$ are removed from $H_{R,h}$, four congruent
``quadrants'' remain. One of the quadrants contains portions of
the positive rays of $X$ and $Z$ in its boundary. Another quadrant
contains portions of the negative rays of $X$ and $Z$ in its
boundary.  We let $Q_{R,h}$ be the union of those two quadrants
(See Figure~\ref{figure:sixpics}). The rotation $\rho_Y$
interchanges the two quadrants, so $Q_{R,h}$ is
$\rho_Y$-invariant. The boundary $\Gamma(R,h)$ of $Q_{R,h}$ is a
piecewise smooth curve that is embedded except at the origin, a
double point.

If in $H_{R,h}$ we replace $Q_{R,h}$ by a connected,
$\rho_Y$-invariant minimal surface $D$ with the same boundary,
then to insure $\rho_Z$ symmetry we must also replace
$\rho_Z(Q_{R,h})$ by the corresponding surface $\rho_Z(D)$.
Because $D$ is $\rho_Y$-invariant,
$\rho_X(D)=\rho_Z(\rho_Y(D))=\rho_Z(D)$.  Thus our candidate
surface is
 \begin{equation}\label{e:T}
    M=M_{R,h}=\overline{D}\cup\rho_Z(\overline{D})
     =\overline{D}\cup\rho_X(\overline{D}).
 \end{equation}

 An Euler characteristic calculation shows that
$M_{R,h}$ has genus one if and only if $D$ is a disk. For instance
if $D$ is a disk, then we can use the nine corners and twelve
edges of the quadrants of $H_{R,h}$ together with the two faces
$D$ and $\rho_Z(D)$ to calculate that the Euler characteristic
$\chi$ of $M$ is $9-12+2=-1$. Since $\partial M = \partial
H_{R,h}$ has one component, its Euler characteristic and genus $g$
are related by $\chi=2-2g-1$. Thus the genus is one, as desired.
(See Figure~\ref{figure:D and M}.)

Note that if $D$ is embedded on one side of $H$, then $\rho_Z(D)$
will lie on the other side of $H$ and thus $M$ will be embedded.
By Schwarz reflection (see~\eqref{e:T}), the interior of $M$ will
be smoothly embedded along $X$ and $Z$ except possibly at the
origin. Since an embedded minimal surface of finite topology in
$\RR^3$ cannot have an isolated interior singularity (\cite{ni2},
\S363), in fact the interior of $M$ is smoothly embedded
everywhere.

  \begfig
 %\hspace{0.3in}
\vspace{.2in} \centerline {
\includegraphics[width=1.5in]{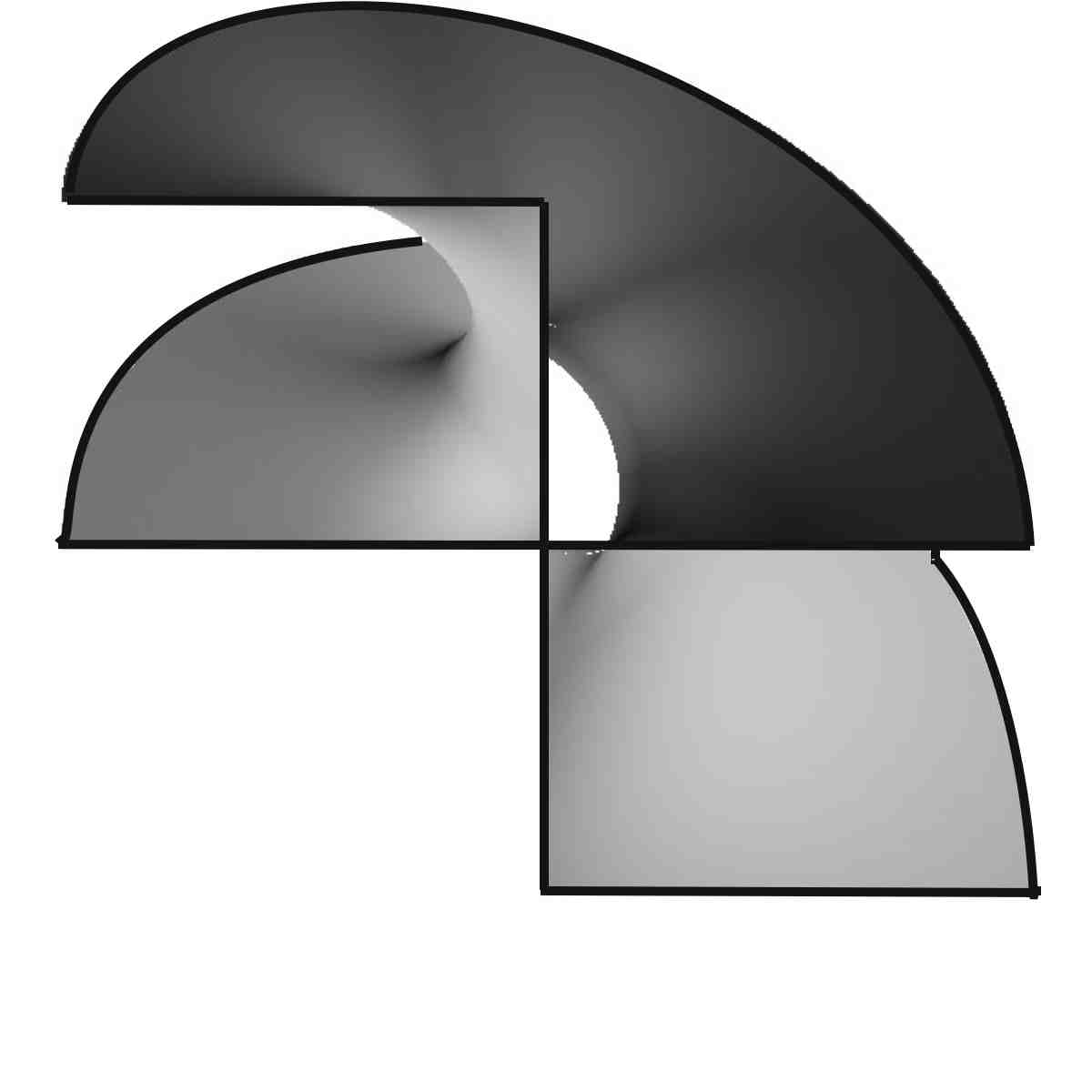}
\includegraphics[width=1.5in]{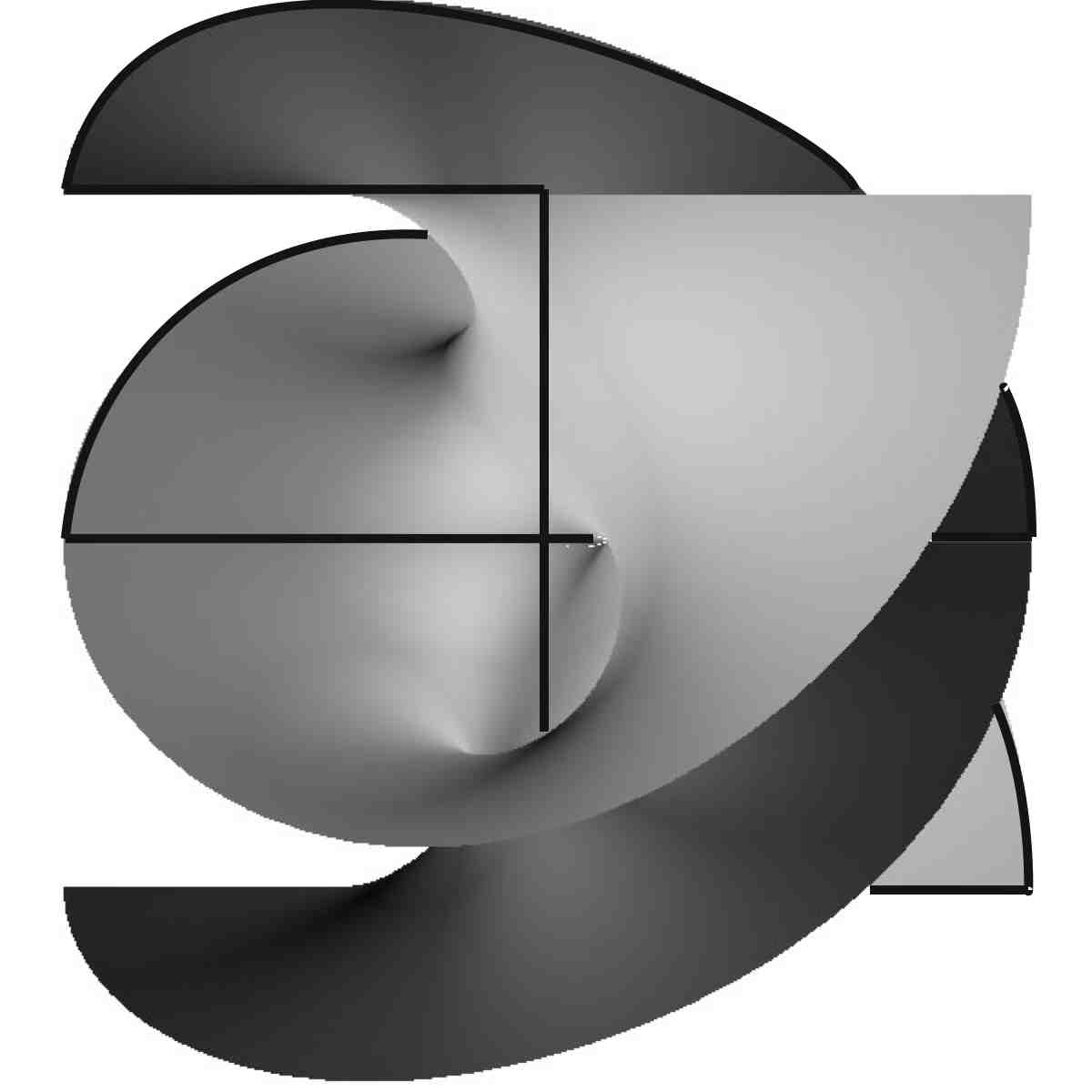}
} \vspace{-.5in}
 \hfill
 \begin{center}
 \vspace{0.7cm}
  \caption{} \label{figure:D and M}\vspace{.3cm}
   \parbox{4.0in}{
   Left: An embedded, $\rho_Y$-invariant  minimal disk $D$ on one side of $H$ with
   $\partial D= \Gamma(R,h)$. In this illustration, $h=\pi$.\newline
 Right: The genus-one surface $M=\overline{D}\cup\rho_Y(\overline {D})$.
  }
 \end{center}
 \endfig

We can extend $M_{R,h}$ to get a $\sigma_{2h}$-invariant surface
$$
  N = N_{R,h} = \bigcup_{n\in \ZZ}\sigma_{2hn}(M_{R,h}).
$$
This surface is clearly embedded, and its boundary consists of two
helices.  Since we can also obtain $N$ from $M$ by repeated
Schwarz reflections about the top and bottom edges,  the interior
of $N$ will be smooth (as was the case with $M$.)  Since $N$ has
locally finite topology and since $\partial N$ is smooth and lies
on the boundary of the convex hull of $N$, there are no boundary
singularities (\cite{ni2}, \S366).  That is, $N$ is smooth
everywhere.

To summarize our discussion, we have the following procedure for
generating symmetric genus-one helicoids:

\begin{itemize}
 \item Find a $\rho_Y$-invariant minimal embedded
 disk $D=D_{R,h}$ that has boundary $\Gamma(R,h)$ and that lies on one side of $H$.
 \item Form the corresponding smooth embedded genus-one surface
 $M=M_{R,h}$ and the corresponding smooth, embedded $\sigma_{2h}$-invariant
 surface $N=N_{R,h}$.
 \item Take a limit of $M_{R,h}$ (or, equivalently, of $N_{R,h}$) as $R$ and
 $h$ tend to infinity to get a nonperiodic example.  Fix $h$ and
 take a limit of $N_{R,h}$ as $R\to\infty$ to get a periodic,
 $\sigma_{2h}$-invariant example.
\end{itemize}
Here and throughout the paper, ``disk'' means ``open disk''.
(Since $\Gamma(R,h)$ has a double point, $\overline{D}$ will not
be embedded even if $D$ is embedded.)  To avoid tedious
repetition, we will let ``disk'' mean ``embedded,
$\rho_Y$-invariant disk'' for the remainder of this section.

\stepcounter{theorem}\subsection{What could go
wrong?}\label{subs:wrong} To prove that the procedure we have just
described works, we must address the following questions:
\begin{enumerate}
 \item How do we know that there is a minimal disk $D$ on one side of $H$ with
 $\partial D=\Gamma(R,h)$?
 \item How do we know that we have smooth convergence as $R$ tends to infinity, or as $R$ and $h$
 both tend to infinity?
 \item How do we know that the limit surface has the desired topology?
\end{enumerate}

Concerning question (1), there is no such minimal disk when $R$ is
very small or when $h\le \pi/2$.  Indeed, in those cases one can
prove that $Q_{R,h}$ is the unique minimal variety with boundary
$\Gamma(R,h)$, and it is not a disk, but rather two disks.
Fortunately, we need minimal disks $D$ only for $h>\pi/2$ and $R$
very large.

We will discuss the curvature estimates that address question (2)
later. The key to questions (1) and (3) turns out to be the
following fact:

\begin{proposition}\label{first proposition}
Let $h>\pi/2$.  For all sufficiently large $R$, there exists a
minimal annulus $A$ that has boundary in $Q_{R,h}$ and that lies
in the component $H^+$ of $\RR^3\setminus H$ containing the
positive $y$-axis.
\end{proposition}

\begin{proof}[Sketch of proof]
(See Proposition~ \ref{prop:aip} for details.) Note that when $R$
is very large, the region $H^+$ near the point $(0,R/2,0)$
resembles a slab between two horizontal planes. Consequently, the
intersection of $H^+$ with a suitable catenoid centered at
$(0,R/2,0)$ has one component that is an annulus $A$ with boundary
in $H$.  The condition $h>\pi/2$ is precisely the condition that
guarantees (provided $R$ is sufficiently large) that $\partial A$
lies in $Q_{R,h}$. (Figure~\ref{figure:barrier} shows the case
$h=\pi$.)
\end{proof}

We can now get a minimal disk in $H^+$ by minimizing area among
disks $D$ that have boundary $\Gamma(R,h)$ and that lie in the
unbounded component of $H^+\setminus A$.

(Since the standard theorems about minimal disks assume embedded
boundary, we should first approximate $\Gamma(R,h)$ by embedded
smooth curves in $H$, but in this outline we will ignore such
technicalities.)

However, that least-area disk $D^*$ is not the one we want.
Consider for example the case $h=\pi$.    By construction, the
disk $D^*$ is disjoint from some catenoid passing through the
middle of $\Gamma(R,h)$, as shown in Figure~\ref{figure:barrier}.
By translating the catenoid around horizontally, we see from the
maximum principle that $D^*$ is forced to lie close to its
boundary.  Consequently, if we take a limit of such $D^*$ as $R\to
\infty$, then we are left just with two flat strips in the
$xz$-plane.  The corresponding complete embedded surface
(generated by Schwarz reflections from those strips) is the
$xz$-plane, not the periodic genus-one helicoid we want.

  \begfig
 \hspace{0.3in}
      \vspace{.2in} \centerline {
       \includegraphics[width=3.0in]{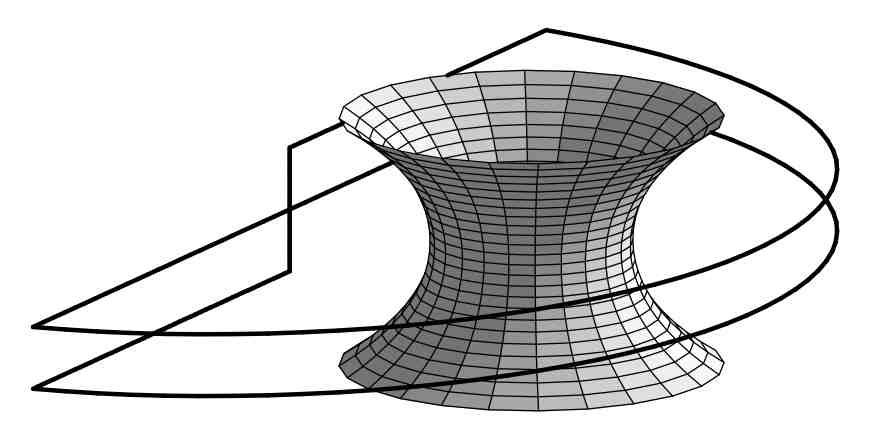}
                                              }
 \hfill
 \begin{center}
 % \vspace{0.5cm}
  \caption{} \label{figure:barrier} \vspace{0.3cm}
        \parbox{4.0in}{
        The curve $\Gamma=\Gamma(R,\pi)$  and a catenoidal barrier. A minimal
         disk that has boundary $\Gamma$ and that avoids the barrier must lie close to $\Gamma$,
         because if we translate the barrier horizontally, it must touch $\Gamma$ before it touches the
         interior of the disk.
         }
 \end{center}
 \endfig

The problem is a loss of topology. The closure of $D^*$ is not
simply connected: it has a closed geodesic starting and ending at
the origin.  But as $R\to\infty$, the length of that geodesic
tends to infinity.  Thus the limiting geodesic is not closed, and
the limiting surface is simply connected.

Similar reasoning shows that in general our procedure is doomed to
yield a simply connected limit surface unless we use minimal disks
$D$ with $\partial D=\Gamma(R,h)$ that have the following {\em
annular intersection property}: $D$ intersects every minimal
annulus $A$ in $H^+$ with $\partial A\subset Q_{R,h}$.

For minimal disks $D$ that have the annular intersection property,
we prove a uniform bound on the length of the geodesic starting
and ending at the origin. This implies that a limit of such disks
will generate a surface with the desired topology. (In particular,
the surface will contain a closed geodesic and therefore will have
nontrivial fundamental group.)

Thus our scheme for producing genus-one helicoids works if and
only if we use disks that have the annular intersection property.
Fortunately such disks do exist:

\begin{proposition}
Let $h>\pi/2$.  For all sufficiently large $R$, the curve
$\Gamma(R,h)$ bounds a $\rho_Y$-invariant minimal embedded disk in
$H^+$ that has the annular intersection property.
\end{proposition}

\begin{proof}[Sketch of proof] (See Section~\ref{section:DiskConstruction} for details.)
By Proposition~\ref{first proposition}, the set $\FF$ of minimal
disks in $H^+$ bounded by $\Gamma(R,h)$ is nonempty.  For
simplicity, let us assume that that $\FF$ is a finite set, and
that each disk in $\FF$ is strictly stable or strictly unstable.
Choose a disk $D\in \FF$ that is closest to $Q_{R,h}$ in the sense
that no other disk in $\FF$ lies between $D$ and $Q_{R,h}$. We
will show that $D$ must have the annular intersection property.

Note that $Q_{R,h}$ is strictly stable: each of its two components
lies in a half-helicoid (i.e., one of the components of
$H\setminus Z$), and each half-helicoid is stable because its
Gauss-map image lies in a hemisphere.

If $D$ were strictly stable, then (by a general minimax or
mountain pass lemma) there would be an unstable disk between $D$
and $Q_{R,h}$, contradicting the choice of $D$. Thus $D$ is
strictly unstable.

Now consider a minimal annulus $A$ in $H^+$ with $\partial A
\subset Q_{R,h}$. If $A$ were disjoint from $D$, we could minimize
area among all disks that have boundary $\Gamma(R,h)$ and that lie
in the region of $H^+$ between $D$ and $A$.  By the instability of
$D$, the result would be a minimal disk $D'$ lying strictly
between $D$ and $Q_{R,h}$, contradicting the choice of $D$.  Thus
$D$ intersects every such annulus.  That is, $D$ has the annular
intersection property.
\end{proof}

We end this outline by saying a word about the curvature estimates
that guarantee smooth convergence in question~(2) of
Section~\ref{subs:wrong} above. The points of a disk $D\subset
\RR^3$ with vertical tangent planes are of course the critical
points of nonzero linear functions of the form $f(x,y,z)=ax+by$.
Morse theory lets us deduce facts about the set of such critical
points from knowledge of the boundary.  In that way, in
Section~\ref{section:Vertical}, we control the set of points in
$D$ with vertical tangent planes. That control in turn lets us
deduce curvature estimates in Section~\ref{section:uniform_est}.

\section{Disks with the annular intersection
property}\label{section:DiskConstruction}

\stepcounter{theorem}\subsection{}\label{s:defs} Let  $H_{R,h}$ be
the intersection of the helicoid $H$ with a solid right-circular
cylinder centered at the origin with axis $Z$, radius $R$, and
height $2h$.  Thus
 $$
   H_{R,h} = F ( (-R,R)\times (-h,h) )
 $$
where
 $$
   F(u,v) = (u\cos v, u\sin v, v).
 $$

The $x$- and $z$-axes divide $H_{R,h}$ into four congruent
``quadrants''.  Let $Q_{R,h}$ denote the union of the first and
third quadrants:
 $$
  Q_{R,h} = F( (0,R)\times(0,h)) \cup F( (-R,0)\times (-h,0)).
 $$
Note that $Q_{R,h}$ consists of two pieces that have a common
corner at the origin and that are related by the $180^\circ$
rotation $\rho_Y$ about the $y$-axis, $Y$. (See
Figure~\ref{figure:sixpics}.)

We let $\Gamma(R,h)$ be the boundary of $Q_{R,h}$.  We will regard
$\Gamma(R,h)$ as a piecewise smooth curve that is embedded except
at the origin, a double point.  Note that $\Gamma(R,h)$ consists
of line segments together with the two helical arcs
\begin{equation}
\label{e:hel arcs, first mention}
 \text{
     $\{ F(R,v): 0\le v \le h\}$
     and
     $\{ F(-R,v): -h\le v\le 0\}$.}
\end{equation}

Recall that $H^+$ is the component of $\RR^3\setminus H$ that
contains the positive $y$-axis. Our goal in this section is to
prove existence of a minimal disk in $H^+$ with boundary
$\Gamma(R,h)$ and with the following ``annular intersection
property'':

\begin{definition}\label{def:annular intersection property}
 If $\Gamma$ is closed curve in
$H$, let $\Aa(\Gamma)$ be the set of minimal embedded annuli $A$
in $H^+$ such that $\partial A$ is smooth and is contained in the
union of the bounded components of $H\setminus\Gamma$. A minimal
surface with boundary $\Gamma$ that intersects  every annulus in
$\Aa(\Gamma)$ is said to have the {\em annular intersection
property}.
\end{definition}

Of course the annular intersection property is vacuous if
$A(\Gamma)$ is empty.  However, $\Aa(\Gamma(R,h))$ is nonempty for
suitable $R$ and $h$:

\begin{proposition}\label{prop:aip}
For every $\eta>\pi/2$, there is an $R_\eta<\infty$ such that
$\Aa(\Gamma(R,h))$ is nonempty provided $R\ge R_\eta$ and $h\ge
\eta$.
\end{proposition}

\begin{proof}
Note that for $h\ge \eta$,
$$
   Q(R,\eta)\subset Q(R,h)
$$
and so $\Aa(\Gamma(R,\eta))\subset \Aa(\Gamma(R,h))$.  Thus it
suffices to prove that $\Aa(R,\eta)$ is nonempty for all
sufficiently large $R$.  For the same reason, it suffices to
consider $\eta$ with $\pi/2< \eta\le \pi$.

Translate $H^+$ and $Q_{R,\eta}$ by $(0,-R/2,0)$ to obtain
$(H^+)'_R$ and $Q'_{R,\eta}$.  Note that as $R\to \infty$,
$(H^+)'_R$ converges to a limit $(H^+)'$ consisting of horizontal
slabs, one of which is the slab $|z|\le \pi/2$.  Also,
$Q'_{R,\eta}$ converges smoothly to the union $Q'$ of the planes
$z=\pi/2$ and $z=-\pi/2$.

Let $C$ be any catenoid with a vertical axis of symmetry.  Note
that $C$ intersects the planes $z=\pm\, \pi/2$ transversely in a
pair of circles that bound an annular portion of $C$ in $(H^+)'$.
Hence for all sufficiently large $R$, the catenoid $C$ intersects
$Q'_{R,\eta}$ in a pair of curves that bound an annular component
$A_R$ of $C$ in $(H^+)_R$.  Translating $A_R$ by $(0,R/2,0)$
produces a minimal annulus in $H^+$ with boundary in $Q_{R,\eta}$.
\end{proof}

\stepcounter{theorem}\subsection{The existence result for smooth
simple closed curves $\Gamma$}

The closed curve $\Gamma(R,h)$ is neither simple nor smooth. We
first prove the result we desire for smooth simple closed curves
$\Gamma$ that approximate $\Gamma(R,h)$.  We then prove the
estimates (Lemma~\ref{first curvature estimate}) that allow us
conclude the desired result (Theorem~\ref{existence3}) for
$\Gamma(R,h)$.

\begin{theorem}\label{minmax}
Suppose $\Gamma$ is a smooth simple closed curve in $H$ such that
\begin{enumerate}
 \item the region $D_\Gamma$ in $H$ bounded by $\Gamma$ is strictly
 stable, and
 \item $\Aa(\Gamma)$ is nonempty.
\end{enumerate}
Then $\Gamma$ bounds a weakly unstable minimal embedded disk $D$
in $H^+$ with the annular intersection property.

If $\Gamma$ is $\rho _Y$-invariant, then we may require that $D$
also be $\rho _Y$-invariant.
\end{theorem}

\begin{proof}
Consider first the case that $\Gamma$ is noncritical in the
following sense: $0$ is not an eigenvalue of the Jacobi operator
of any smooth embedded minimal disk in $H^+$ with boundary
$\Gamma$.

Let $\FF=\FF(\Gamma)$ be the set of minimal embedded disks in
$H^+$ bounded by $\Gamma$.  Then $\FF\cup\{D_{\Gamma}\}$ is
compact by standard curvature estimates (see for
example~\cite{wh5} or Lemma~\ref{first curvature estimate} below.)
Moreover, since $D_\Gamma$ is strictly stable and since $\Gamma$
is noncritical, $\FF$ is in fact finite.

Let $A$ be an annulus in $\Aa(\Gamma)$. We claim that $\Gamma$
bounds an embedded disk in $H^+$ disjoint from $A$.  To see this,
let $\BB$ be a large ball in $\RR^3$ centered at the origin with
$\Gamma$ in its interior.  Note that $H\cap \BB$ and $H^+\cap
(\partial \BB)$ are topologically disks with a common boundary, so
their union $S$ is topologically a sphere.  Thus $\Gamma$ divides
$S$ into two regions, $D_{\Gamma}$ and $S\setminus D_{\Gamma}$,
each of which is topologically a disk. In particular, $S\setminus
D_{\Gamma}$ is a piecewise smooth embedded disk in
$\overline{H^+}$ that has boundary $\Gamma$ and that is disjoint
from $A$.  By perturbing slightly, we get a smooth embedded disk
in $H^+\setminus A$ with boundary $\Gamma$.

Now minimize area among all disks in $H^+\setminus A$ with
boundary $\Gamma$. The minimum exists and is smoothly embedded by
a theorem of Meeks and Yau (see Theorem~\ref{Meeks-Yau} below),
and thus it is a disk in the family $\FF$. By finiteness of $\FF$,
there is a disk $D\in \FF$ that is closest to $D_{\Gamma}$ in the
sense that no other disk in $\FF$ lies between $D$ and
$D_{\Gamma}$.

Now if $D$ were stable, it would be strictly stable by
noncriticality of $\Gamma$.  But then by a standard minimax
principle (see for example Theorem \ref{minimax theorem in
appendix} in the appendix), there would be an unstable minimal
embedded disk between $D$ and $D_{\Gamma}$, contradicting the
choice of $D$. Thus $D$ is strictly unstable.

It remains only to show that $D$ must intersect every annulus in
$\Aa(\Gamma)$.  Suppose on the contrary that ${\mathcal
A}(\Gamma)$ contains an annulus $A$ disjoint from $D$.  Let $D'$
be the least-area embedded disk bounded by $\Gamma$ in the closure
of the region of $H^+$ between $A$ and $D$. (The disk exists and
is smoothly embedded by the Meeks-Yau theorem.) Since $D_{\Gamma}$
does not lie in the closure of that region, $D'\ne D_{\Gamma}$.
Since $D$ is strictly unstable, $D'\ne D$. Thus $D'$ lies between
$D$ and $D_{\Gamma}$, contradicting the choice of $D$. Hence $D$
intersects $A$ as claimed.

This completes the proof in the case of noncritical $\Gamma$. In
fact, noncritical $\Gamma$ are generic (see Theorem~\ref{degree
theory theorem} in the appendix). Thus we can find a nested
sequence of noncritical $\Gamma_i\subset H$ converging smoothly to
$\Gamma$ from the outside. Let $D_i$ be an unstable embedded disk
in $H^+$ bounded by $\Gamma_i$ and intersecting all the $A\in
\Aa(\Gamma_i)$. Note that $\Aa(\Gamma)\subset {\mathcal
A}(\Gamma_i)$ since $\Gamma_i$ encloses $\Gamma$ in $H$. Thus
$D_i$ intersects every $A\in \Aa(\Gamma)$. The curvatures of the
$D_i$ are bounded by standard estimates (see \cite{wh5} or
Lemma~\ref{first curvature estimate} below.) Thus a subsequence
converges smoothly to a limit disk $D$ in $\overline{H^+}$. Now
$D$ is weakly unstable (since it is the smooth limit of unstable
disks), so $D\ne D_\Gamma$.  Hence (by the strong maximum
principle) $D$ cannot touch $H$, so $D$ lies in $H^+$. Also, $D$
intersects every $A\in \Aa(\Gamma)$ since each $D_i$ does. This
completes the proof for arbitrary $\Gamma$.

Finally, the proof in case of $\rho_Y$-invariance is exactly the
same, except that we work throughout with $\rho_Y$-invariant
curves and disks, and with the restriction of the Jacobi operator
to the space of $\rho_Y$-equivariant vectorfields.  The minimal
annulus $A$ should be replaced by $A\cup \rho_Y(A)$. Where the
proof uses a least-area disk $\Delta$, that disk turns out to be
$\rho_Y$-invariant.  For if $\Delta$ were not $\rho_Y$-invariant,
then $\rho_Y(\Delta)$ would be a second least-area disk with the
same boundary.  By the Meeks-Yau theorem, $\Delta$ and
$\rho_Y(\Delta)$ would be disjoint.  But by the $\rho_Y$ symmetry,
the volume of the region between $D_{\Gamma}$ and $\Delta$ must
equal the volume of the region between $D_{\Gamma}$ and
$\rho_Y(\Delta)$, so $\Delta$ and $\rho_Y(\Delta)$ cannot be
disjoint.  The contradiction proves that $\Delta$ is indeed
$\rho_Y$-invariant.
\end{proof}

For the reader's convenience, we state the theorem of Meeks and
Yau that was used in the preceding proof:

\begin{theorem}[Meeks-Yau]\cite{my1,my2}\label{Meeks-Yau}
Let $\Omega\subset \RR^3$ be a mean convex domain with piecewise
smooth boundary, and let $\Gamma$ be a smooth curve in $\partial
\Omega$ that bounds a disk in $\Omega$.  Then $\Gamma$ bounds a
least-area disk $D$ in $\overline{\Omega}$.  Such a disk must be
smooth and embedded, and it must be contained either in $\partial
\Omega$ or else in $\Omega$.  Furthermore, any two such disks must
be disjoint.
\end{theorem}

Meeks and Yau prove this theorem for convex domains in Theorem 6
of \cite{my1}. They extend the result to other domains in
Section~1 of \cite{my2}.

\stepcounter{theorem}\subsection{Approximation results and uniform
curvature estimates}\label{approx}

Note that a half-helicoid (such as either component of $H\setminus
Z$) is stable since its image under the Gauss map is contained in
a hemisphere.  Consequently any bounded domain in a half-helicoid
is strictly stable.  In particular, the union $Q=Q_{R,h}$ of the
bounded components of $H\setminus \Gamma(R,h)$ is strictly stable.
In the following lemma, we show that $Q$ can be fattened to get a
strictly stable domain $\Omega\subset H$ bounded by a smooth
simple closed curve.

\begin{lemma}\label{stable}
For any positive numbers $R$ and $h$, there is a simply connected,
$\rho_Y$-invariant domain $\Omega\subset H$ that contains
    $\overline{Q_{R,h}}$,
    is strictly stable, and
    is bounded by a smooth embedded  curve.
\end{lemma}

\begin{proof}
Let $\Omega_i$ be the set of points in $H$ at distance $\le 1/i$
from the union $Q=Q_{R,h}$ of the two bounded components of
$H\setminus \Gamma$.  Let $\lambda_i$ and $\lambda$ be the lowest
eigenvalues of the Jacobi operator on $\Omega_i$ and on $Q$,
respectively. Let $u_i$ be the first eigenfuction on $\Omega_i$,
normalized so that $\max u_i=1$. Recall that $u_i>0$ at every
interior point. Note that
$$
   \lambda_1 < \lambda_2 < \lambda_3 < \dots
$$
and that $\lambda_i<\lambda$ for all $i$, so that
\begin{equation}\label{e:eigenbound}
  \lambda_\infty \le \lambda
\end{equation}
where
$$
  \lambda_\infty=\lim_i\lambda_i.
$$
By the Schauder estimates, we may assume that the $u_i$ converge
to a limit $u$, the convergence being smooth away from the corners
of $\Gamma$.

(If the smooth convergence is not clear, note that the $u_i$ are
uniformly bounded since they are normalized to have maximum value
$1$.  Thus the Schauder estimates give uniform local
$C^{2,\alpha}$ bounds as $i\to\infty$ away from the corners of
$\Gamma$. Those bounds imply $C^2$ convergence away from the
corners for a subsequence of the $u_i$. Likewise, higher-order
Schauder estimates imply $C^k$ convergence for each $k$.)

\begin{claim} The convergence $u_i\to u$ is uniform up to and
including the boundary.
\end{claim}

\begin{proof}[Proof of claim]
It suffices to show the following: if $p_i\in \Omega_i$ converges
to $p\in \Gamma$, then $u_i(p_i)$ converges to $0$.  To see that
it does, fix a $p\in \Gamma$ and let $a$ be the supremum of
 $$
   \limsup u_i(p_i)
 $$
among all sequences $p_i\in \Omega_i$ with $p_i\to p$.  Note that
the supremum is attained by some sequence $p_i$. By passing to a
subsequence, we may assume that
 $$
   \lim u_i(p_i) = a.
 $$
Our goal is to show that $a=0$.

Translate $\Omega_i$ by $-p_i$ and dilate by
  $$
    \mu_i := \frac1{\dist(p_i,\partial \Omega_i)}
  $$
to get a surface $\Omega_i'$.  Note that
\begin{equation}\label{e:mu}
    \mu_i\to \infty.
\end{equation}
By passing to a subsequence, we can assume that the $\Omega_i'$
converge to a plane domain $\Omega'$ with
$$
  \dist(0,\partial \Omega') = 1
$$
so that in particular $\partial \Omega'$ is not empty.  Note that
$\partial \Omega'$ is either smooth or piecewise smooth.  (Indeed,
$\partial \Omega'$ must be a straight line, or two rays together
with a quarter circle joining their endpoints, or a right angle,
or the union of two disjoint right angles.  Here ``right angle''
means ``union of two orthogonal rays with a common endpoint''.)

Let $u_i'$ be the function on $\Omega_i'$ corresponding to $u_i$.
Then $u_i'$ is a Jacobi eigenfuction with eigenvalue
 $$
   \lambda_i' = \frac{\lambda_i}{\mu_i^2}.
 $$
Note that $\lambda_i'\to 0$ (because by~\eqref{e:mu} the $\mu_i$
tend to infinity and by~ \eqref{e:eigenbound} the $\lambda_i$'s
are bounded.) Thus by passing to a subsequence we may suppose that
the $u_i'$ converge smoothly away from the corners of $\partial
\Omega'$ to a Jacobi eigenfunction
 $$
   u': \Omega' \to \RR
 $$
with eigenvalue $0$.  Since $\Omega'$ is planar, $u'$ is in fact a
harmonic function.  Note that
$$
  \max u' = u'(0) = a.
$$
Thus by the maximum principle for harmonic functions, $u'$ is
constant. But $u'$ vanishes on the smooth portions of $\partial
\Omega'$, so $a$ must be $0$. This completes the proof of the
claim.
\end{proof}

We now resume the proof of the lemma.  By the uniform convergence
$u_i\to u$, $u$ is nonzero (its maximum value is $1$) and it
vanishes at the boundary.  By the smooth convergence on the
interior, it is a Jacobi eigenfunction on $Q$ with eigenvalue
$\lambda_\infty$.   Since $Q$ is strictly stable (see the
discussion immediately preceding the lemma), $\lambda_\infty>0$.

Thus there is an $n$ (any sufficiently large $n$ will do) for
which $\Omega_n$ is strictly stable.  Now let $\Omega$ be any
$\rho_Y$-invariant domain in $H$ such that
$$
  \overline{Q} \subset \Omega \subset \Omega_n
$$
and such that $\partial\Omega$ is a smooth simple closed curve.
\end{proof}

The next lemma follows from the estimates in \cite{wh5}, but since
the proof is short, we include it for the reader's convenience.

\begin{lemma}\label{first curvature estimate}
Let $D_i$ be a sequence of
embedded minimal disks in $H^+$ with embedded, piecewise smooth
boundary curves $\Gamma_i$ in $H$. Suppose that the $\Gamma_i$
converge to a curve $\Gamma$ that is smooth and embedded except at
a finite set $S$, and suppose that the convergence
$\Gamma_i\to\Gamma$ is smooth on compact subsets of
$\RR^3\setminus S$.  Suppose also that the total curvatures of the
$\Gamma_i$ are uniformly bounded.

Then the principal curvatures of the $D_i$ are uniformly bounded
as $i\to\infty$ on compact subsets of $\RR^3\setminus S$.
\end{lemma}

\begin{proof}  In this proof, it will be convenient to let $D_i$
denote the disk together with its boundary. Suppose the lemma
fails.  Then (by passing to a subsequence) there is a compact set
$K$ disjoint from $S$ and a sequence $p_i\in K\cap \overline{D}_i$
such that
 $$
      B(\overline{D}_i,p_i)\to \infty.
 $$
Here $B(\overline{D}_i,p_i)$ is the norm of the second fundamental
form of $\overline{D}_i$ at $p_i$.  By enlarging $K$ slightly, we
can assume that
 \begin{equation}\label{e:blowup}
     B(\overline{D}_i,p_i)\, \dist(p_i,\partial K) \to \infty.
 \end{equation}
Fixing this $K$, we may rechoose $p_i\in K\cap \overline{D}_i$ to
maximize the left hand side of~\eqref{e:blowup}.

Now translate $D_i$, $K$, and $H$ by $-p_i$ and dilate by
$$
  \mu_i := B(\overline{D}_i, p_i)
$$
to get $D_i'$, $K_i'$, and $H_i'$ such that
\begin{equation}\label{e:normalization}
  B(\overline{D}_i{}',0) = 1
\end{equation}
and such that
\begin{equation*}
    \max_{p\, \in K_i'\cap \overline{D}_i{}'}
    B(\overline{D}_i{}', p)\, \dist(p,\partial K_i')
    =
    \dist(0, \partial K_i') \to \infty.
\end{equation*}
From this we see that the $\partial K_i'$ are moving off to
infinity and that the curvatures of the $D_i'$ are uniformly
bounded as $i\to\infty$ on compact subsets of $\RR^3$. Note that
$\partial D_i'$ converges either to the empty set (i.e., it moves
off to infinity) or to a straight line, the convergence in the
latter case being smooth on compact subsets of $\RR^3$. Thus we
may assume that the $D_i'$ converge smoothly to a limit minimal
surface $D'$.  After passing to a subsequence, the $H_i'$ will
converge smoothly to a limit $H'$ that is either a plane or the
empty set.

Now $D'$ is simply connected and embedded, and
 \begin{equation}\label{e:nonflat}
   B(D',0)=1
 \end{equation}
by~\eqref{e:normalization}.  By the Gauss-Bonnet theorem, the
total curvatures of the $D_i$ are uniformly bounded and thus $D'$
has finite total curvature.

If $D'$ has no boundary, then it is a complete, embedded, simply
connected minimal surface of finite total curvature.  But the only
such surface is a plane, contradicting~\eqref{e:nonflat}.

Thus $D'$ has nonempty boundary, namely a line.  Note that the
line lies in a plane, namely $H'$, and that $D'$ lies in a closed
halfspace bounded by $H'$.  Thus extending $D'$ by Schwarz
reflection produces a complete, simply connected, embedded minimal
surface of finite total curvature.  Again, the only such surface
is a plane, contradicting~\eqref{e:nonflat}.
\end{proof}

\stepcounter{theorem}\subsection{The existence result for
$\Gamma(R,h)$} \label{D_Rh} Recall that ${\mathcal  A}(\Gamma)$ is
the set of minimal annuli $A$ in $H^+$ with $\partial A$ contained
in the union of the bounded components of $H\setminus \Gamma$, and
that $\Gamma(R,h)$ is the piecewise smooth curve defined in
Section~\ref{s:defs}.

\begin{theorem}\label{existence3}
Let $\Gamma=\Gamma(R,h)$ and suppose $\Aa(\Gamma)$ is not empty.
Then $\Gamma$ bounds a $\rho_Y$-invariant minimal disk $D$ in
$H^+$ such that
 \begin{enumerate}
 \item $\overline{D}$ is smoothly embedded except at the corners
 of $\Gamma$, and
 \item $D$ has the annular intersection property
 (\ref{def:annular intersection property}).
 \end{enumerate}
\end{theorem}

\begin{proof}
By Lemma~\ref{stable}, $H$ contains a $\rho_Y$-invariant, simply
connected, strictly stable domain $\Omega$ containing $\Gamma$
such that $\partial\Omega$ is a smooth embedded curve.

Now let $\Gamma_i$ be a sequence of smooth $\rho_Y$-invariant
simple closed curves in $\Omega$ such that
 \begin{enumerate}
 \item $\Gamma_i$ encloses $\Gamma$,
 \item $\Gamma_i$ converges to $\Gamma$ as $i\to \infty$, the
   convergence being smooth away from the corners of $\Gamma$, and
 \item the total curvatures of the $\Gamma_i$ are uniformly
 bounded.
 \end{enumerate}

Since the region of $H$ bounded by $\Gamma_i$ is a subset of
$\Omega$, it must be strictly stable. Since $\Gamma_i$ encloses
$\Gamma$, $\Aa(\Gamma_i)$ contains $\Aa(\Gamma)$ and is therefore
nonempty.

Thus by Theorem~\ref{minmax}, $\Gamma_i$ bounds a
$\rho_Y$-invariant minimal embedded disk $D_i$ in $H^+$ that
intersects every annulus $A$ in $\Aa(\Gamma_i)$.  In particular,
$D_i$ intersects every annulus $A$ in $\Aa(\Gamma)$.

By Lemma~\ref{first curvature estimate}, the curvatures of the
$D_i$ are uniformly bounded away from the corners of $\Gamma$.
Thus by passing to a subsequence we may assume that the $D_i$
converge smoothly (away from those corners) to an embedded minimal
surface $D$ in the closure of $H^+$ with $\partial D =\Gamma$.

We claim that $D$ is not contained in $H$.  To see that this is
the case, note that the strict stability of $\Omega$ implies that
$\Omega$ is contained in an open set $W$ of $\RR^3$ with the
following property: any minimal surface with boundary in $\Omega$
either is entirely contained in $\Omega$ or else contains points
in $W^c$. (See Corollary~\ref{cor:appendix} in the appendix.)
Since $\partial D_i$ is contained in $\Omega$ and $D_i$ is not
contained in $\Omega$, we see that $D_i\cap W^c$ is not empty.
Thus $D\cap W^c$ is not empty, which implies that $D$ is not
contained in $H$, as claimed.

Since the $D_i$ are simply connected, each component of $D$ is
simply connected. We claim that $D$ consists of only one connected
component. To see that it does, note that $\Gamma$ may be regarded
as the union of two simple closed curves $C'$ and $C''$ that
intersect at the origin. Suppose that $D$ is not connected. Then
it must consist of two components $D'$ and $D''$ with boundary
curves $C'$ and $C''$, respectively. Note that $D''=\rho_Y(D')$ by
the $\rho_Y$ symmetry of $D$. Since $D$ does not lie in $H$,
neither $D'$ nor $D''$ can lie in $H$. Note that $C'$ lies in the
closure $\Sigma'$ of one of the two connected component of
$H\setminus Z$: it consists of a segment $I$ of the $z$-axis, two
horizonal segments, and a helical curve. Now rotate the other
half-helicoid $\Sigma''=H\setminus\Sigma'$ about the $z$-axis in
$H^+$ until it first either touches $D'$ on the interior or
becomes tangent to $D'$ along the $z$-axis.  It cannot touch at an
interior point by the maximum principle.  It cannot touch at the
endpoints of $I$ since $D'$ is tangent to $H$ there. Thus the
first point of contact is a point of tangency at a point inside
the segment $I$. But that violates the boundary maximum principle.
This contradiction proves that $D$ is connected.

Since $D$ lies in the closure of $H^+$ but does not lie in $H$, by
the strict maximum principle $D$ cannot touch $H$ at any interior
point.  That is, $D$ is contained in $H^+$.

Finally, $D$ intersects every $A\in \Aa(\Gamma)$ since each $D_i$
does.
\end{proof}

\section{Vertical tangent planes}
\label{section:Vertical}

Throughout this section, we will assume that $R$ and $h$ are fixed
with $h>\pi/2$, and that $D\subset H^+$ is  an embedded,
$\rho_Y$-invariant minimal disk with boundary
  $\partial D = \Gamma=\Gamma(R,h)$,
the curve specified in Section~\ref{s:defs}.  As explained in
Section~\ref{subsection:replacing},
$$
   M := \overline{D}\cup \rho_Z(\overline{D})
      = \overline{D}\cup \rho_X(\overline{D})
$$
is a smooth, embedded genus-one minimal surface, and
$$
  N := \bigcup_{n\in \ZZ} \sigma_{2nh}(M)
$$
is a smoothly embedded, $\rho_{2h}$-invariant minimal surface
whose boundary is a pair of helices.

We will prove in Proposition~\ref{vertical thm 1} that there are
at most two interior points of $D$ at which the tangent plane is
parallel to a given vertical plane.  Furthermore, such points must
be close (within distance $2\pi$) to the $xy$-plane.    We will
give similar bounds on the number of times a vertical plane $V$
can be tangent to $D$ along $Z$ (Theorem~\ref{VtangentZ}). We
combine these results to get a local upper bound
(Theorem~\ref{Vestimate}) on the number of points in $N$ where the
tangent plane is parallel to a given vertical plane $V$.

The theorems of this section, all of which control the set of
points with vertical tangent planes, will be used in
Section~\ref{section:uniform_est} to get the curvature estimates
needed when we let $R$ and/or $h$ tend to infinity.

A key ingredient in the proof of Proposition~\ref{vertical thm 1}
is a generalization of Rado's theorem, according to which the
tangent plane to a minimal disk in $\RR^3$ must intersect the
boundary in at least four points.  Rado's theorem has the
following consequence: if the restriction of a linear function $f$
to the boundary of a minimal disk has only one local maximum, then
$f$ has no interior critical points in the disk.  The
generalization, due to Schneider (\cite{schneider1} or
\cite{ni2}, {\S 374}), is the following proposition:

\begin{proposition}[Schneider]\label{schneider}
Suppose that $D$ is a minimal disk in $\RR^3$ and that
$f:\RR^3\to\RR$ is a linear function.  If $f\vert \partial D$ has
at most $n$ local maxima with $f>0$, then $f\vert D$ has at most
$(n-1)$ interior critical points, counting multiplicity, with
$f\ge 0$.
\end{proposition}

\stepcounter{theorem}\subsection{A global angle function on $H^+$}
\label{s:angle function}

Let $S$ be the half-strip in the $xz$-plane defined by

\[
 S=\{(x,0,z): x>0,\, -\pi< z < 0\}.
\]
Note that
\[
    H^+
        =
    \bigcup_{_{\theta\in {\RR}}}\sigma_\theta(S),
\]
 where $\sigma_\theta:\RR^3\to\RR^3$ is the screw motion defined by
\[
  \sigma_\theta(r \cos \phi,r\sin \phi,z)=
(r\cos (\phi +\theta) ,r\sin(\phi+\theta),z+\theta).
\]
Thus $H^+$ is foliated by the $\sigma_\theta(S)$, and we may think
of $\theta$ as a globally defined angle function on $H^+$.
Furthermore, $\theta$ extends continuously to
$\overline{H^+}\setminus Z$. Indeed,
 \[
   \theta: \overline{H^+}\setminus Z \to \RR
 \]
is the unique continuous function such that
 \begin{equation*}
   (x,y,z)
    = (\sqrt{x^2+y^2}\cos\theta, \sqrt{x^2+y^2}\sin\theta,z)
 \end{equation*}
and such that
 \begin{equation*}
   \text{$\theta(x,0,0)= 0$ for $x>0$.}
 \end{equation*}
Note also that $\theta(x,0,0)=\pi$ for $x<0$ and that
\begin{equation*}
  \theta(\sigma_s(p))=\theta(p)+s.
\end{equation*}

It will be useful to understand the helical portions of the
boundary curve $\Gamma=\Gamma(R,h)$ (which was defined in
Section~\ref{s:defs}) in terms of the angle function $\theta$.
By~\eqref{e:hel arcs, first mention} in Section~\ref{s:defs},
there are two such helical portions:
\begin{equation}\label{e:top part}
  (R\cos v, R\sin v, v), \qquad 0 \le v \le h
\end{equation}
and
\begin{equation}\label{e:bottom part}
  (-R\cos v, -R\sin v, v), \qquad -h \le v \le 0.
\end{equation}
Now
$$
  (R\cos v, R\sin v, v) = \sigma_v (R,0,0)
$$
so
$$
 \theta(R\cos v, R\sin v, v) = \theta(R,0,0)+v = v.
$$
Similarly
$$
 (-R\cos v, -R\sin v, v) = \sigma_v(-R,0,0)
$$
so
$$
 \theta(-R\cos v,-R\sin v, v)
 =
 \theta(-R,0,0)+v=\pi+v.
$$
Hence we can express the upper and lower helical portions of
$\Gamma$ in terms of the angle coordinate $\theta$ by substituting
$\theta$ for $v$ in~\eqref{e:top part} and by substituting
$\theta-\pi$ for $v$ in~\eqref{e:bottom part}, yielding
\begin{equation}\label{e:upper}
 (R\cos\theta, R\sin\theta, \theta), \qquad 0 \le \theta \le h
\end{equation}
for the upper portion and
\begin{equation}\label{e:lower}
 (R\cos\theta, R\sin\theta, \theta-\pi),
              \qquad -h+\pi \le \theta \le \pi
\end{equation}
for the lower portion.

 The curve $\Gamma$ also contains six
line segments: the four horizontal segments joining the endpoints
of the helical arcs~\eqref{e:upper} and~\eqref{e:lower} to the
$z$-axis, and two segments in the $z$-axis (the two segments given
by $0<z<h$ and $-h<z<0$).  See Figure~\ref{figure:sixpics}.

\stepcounter{theorem}\subsection{Vertical tangent planes at points
in the interior of $D$} \label{interiorvertical}

\begin{proposition}\label{vertical thm 1}
Let $D$ be a $\rho_Y$-invariant minimal embedded disk in $H^+$
with
  $\partial D=\Gamma=\Gamma(R,h)$.
Let
$$
   f(x,y,z) = ax + by
$$
be a nonzero linear function with $a\le 0$.

If $a<0$, then
\begin{enumerate}
 \item\label{item1} $f\vert D$ has at most one critical point $p$ with
 $f(p)\ge 0$.  Such a critical point must satisfy $0<\theta(p)<2\pi$.
 \item\label{item2} $f\vert D$ has at most one critical point $p$ with $f(p)\le
 0$.  Such a critical point must satisfy $-\pi<\theta(p)<\pi$.
\end{enumerate}
If $a=0$, then $f|D$ has exactly one critical point.  That
critical point is the unique point of intersection of $D$ and $Y$,
and it lies on the positive $y$-axis.
\end{proposition}

\begin{remark}\label{remark:asymmetry} There appears to be an asymmetry
between assertions~(\ref{item1}) and~(\ref{item2}), but in fact
they are equivalent (either one implies the other). To see this,
let $g(x,y,z)=ax-y$. If $p$ is a critical point of $f\vert D$ with
$f(p)\le 0$, then $q=\rho_Y(p)$ is a critical point of $g\vert D$
with $g(q)\ge 0$, and $\theta(q)=\pi-\theta(p)$.  The apparent
asymmetry disappears if the theorem is restated using the angle
function
   $\omega(p) := \theta(p) - \pi/2$.
(Note that $\omega(\rho_Y(p))=-\omega(p)$.)
\end{remark}

\begin{proof}
For $k\in \ZZ$, let
$$
  D_k = \{ p \in D: 2\pi k < \theta(p) < 2\pi (k+1)\}.
$$
We will prove that $f\vert D_k$ has
 \begin{enumerate}
 \item[(1a)] no critical points with $f\ge 0$ if $k\ne 0$, and
 \item[(1b)] at most one critical point with $f\ge 0$ if $k=0$.
 \end{enumerate}
These imply assertion~(\ref{item1}), because if $p=(x,y,z)\in
D\setminus\cup_kD_k$, then $\theta(p)$ is an even multiple of
$2\pi$, which implies that $y=0$ and $x>0$ and thus, in the case
$a<0$, that $f(p)<0$. (At the end of the proof, we will also
use~(1a) and~(1b) in the case $a=0$.)

Let us first suppose that $h$ is an integer multiple of $2\pi$:
$$
  h = 2\pi N.
$$
Let $J$ be the halfplane
$$
  J = \{ (x,0,z): x \ge 0\}
$$
and let
$$
  \Gamma_k = (\partial D_k)\setminus J.
$$

Note that if $k\ge N$ or if $k<-N$, then $\Gamma_k=\emptyset$,
which implies by the convex hull property that $D_k=\emptyset$.
Thus from now on we assume that $-N\le k<N$.

By hypothesis, $f\le 0$ on $J$, so any local maxima of $f\vert
\partial D_k$ with $f>0$ must lie on $\Gamma_k$.

If $k>0$, then $\Gamma_k$ consists of the single helical arc:
$$
   \{ (R\cos\theta,R\sin\theta,\theta): 2\pi k <\theta<2\pi (k+1) \}.
$$
Likewise if $k<0$, then $\Gamma_k$ consists of the single helical
arc
$$
   \{ (R\cos\theta,R\sin\theta, \theta - \pi): 2\pi k <\theta<2\pi (k+1) \}.
$$
(These assertions about $\Gamma_k$ follow immediately from
\eqref{e:upper} and \eqref{e:lower} in Section~\ref{s:angle
function}.) Orthogonal projection to the $xy$-plane maps each of
these helical arcs homeomorphically to the circle $x^2+y^2=R^2$
minus the point $(R,0,0)$.   Thus if $k\ne 0$, then
$f\vert\partial D_k$ has exactly one local maximum with $f>0$. By
Schneider's theorem, $f\vert D_k$ has no critical points with
$f\ge 0$. This proves~(1a).

To prove~(1b), note that $\Gamma_0$ has two connected components.
One component is the helical arc
  $$
    \{ (R\cos\theta,R\sin\theta,\theta): 0<\theta<2\pi \}.
  $$
The other connected component of $\Gamma_0$ is the helical arc
   $$
   \{ (R\cos\theta,R\sin\theta,\theta-\pi): 0<\theta<\pi\}
   $$
together with the line segment
   $$
    \{ (t,0,0): 0 > t > -R \}.
   $$
Hence the function $f$ has at most two strictly positive local
maxima on $\partial D_0$. (If this is not clear, consider the
projections to the $xy$-plane as above.) By Schneider's theorem,
$f\vert D_0$ has at most one critical point with $f\ge 0$.  This
proves~(1b).

This completes the proof of~(1a) and~(1b) under the hypothesis
that $h/(2\pi)$ is an integer.  If $2\pi (N-1) < h < 2\pi N$, then
the topmost component of $\Gamma\setminus J$ is not a full turn of
a helix, but rather the helical arc
$$
 \{ (R\cos\theta, R\sin\theta, \theta): 2\pi (N-1) < \theta < h\}
$$
together with the line segment
$$
  \{ (t\cos h, t\sin h, h): 0 < t < R \}
$$
Note that this piecewise smooth arc $C$ still has the following
property: the function $f(x,y,z)=ax+by$ has at most one local
maximum on $C$ with $f>0$. A similar remark applies to the lowest
component of $\Gamma\setminus J$.  The rest of the proof  of~(1a)
and of~(1b) is exactly as before.

As explained above, (1a) and~(1b) imply assertion~(\ref{item1}),
which in turn (by Remark~\ref{remark:asymmetry}) implies
assertion~(\ref{item2}).

Now suppose $a=0$. Without loss of generality we may assume that
$f(x,y,z)=y$. Since $\rho_Y:D\to D$ is an orientation-preserving
involution, it has a unique fixed point $q$.  That is, $q$ is the
unique point of $D\cap Y$.  Since $D\in H^+$, the point $q$ lies
on the positive $y$-axis.  Since $\rho_Y:D\to D$ preserves
orientation, so does $\rho_Y:\Tan_qD\to\Tan_qD$.  Thus $\Tan_qD$
is perpendicular to $Y$, so $q$ is a critical point of $f\vert D$.

By (1a), $q$ is the only critical point of $f\vert \cup_kD_k$ with
$f\ge 0$. In particular, there is no critical point $p$ of $f|D$
for which $\theta(p)$ an odd multiple of $\pi$.  By the $\rho_Y$
symmetry, there is also no critical point $p$ for which
$\theta(p)$ an even multiple of $\pi$.  Thus $q$ is the only
critical point of $f\vert D$ with $f\ge 0$.

Finally, for each $k$ (including $k=0$),
$$
  \Gamma_k \cap \{y<0\}
$$
consists of just one arc (unless it is empty).  (This is because
$y\ge 0$ on one of the two components of $\Gamma_0$.) As before,
$f$ has a unique local minimum on each such arc. Thus by Rado's
theorem (or by Schneider's theorem applied to the function $-f$)
$f\vert D$ has no critical points with $f<0$.
\end{proof}

\begin{corollary}\label{cor:vertical}
Suppose $D\subset H^+$ is a $\rho_Y$-invariant minimal embedded
disk with boundary $\Gamma$. Let $V$ be a plane containing the
$z$-axis.  If $p\in D$ is a point at which $\Tan_pD$ is parallel
to $V$, then $p$ lies within distance $2\pi$ of the plane $z=0$.
There are at most two such points.  If there are two, they lie on
opposite sides of $V$.  If there is such a point on $V$, then it
satisfies:
 \begin{equation}\label{e:PiBound}
  0 < \theta(p) < \pi.
 \end{equation}
\end{corollary}

\begin{proof}  Note that $V$ can be expressed as the zero set of
a nonzero linear function $f(x,y,z)=ax+by$ with $a\le 0$, and that
$p$ is a critical point of $f|D$ if and only if $\Tan_pD$ is
parallel to $V$.  By Proposition~\ref{vertical thm 1}, there are
at most two critical points, and if there are two, they must lie
on opposite sides of $V$.  Also, if $p$ is a critical point, then
$$
   -\pi < \theta(p) < 2\pi
$$
which implies that $p$ lies in the slab $|z|<2\pi$, since, for any
$q=(x,y,z)\in H^+$
$$
     \theta(q)-\pi < z < \theta(q).
$$
Finally, if $p$ is a critical point on $V$, then $f(p)=0$, so
combining conclusions~(1) and~(2) of the theorem
gives~\eqref{e:PiBound}.
\end{proof}

\begin{remark}\label{remark:deformed helical arcs}
Proposition~\ref{vertical thm 1} and Corollary~\ref{cor:vertical},
and indeed all the results in
 Section~\ref{section:Vertical},
are true for a larger class of curves $\Gamma$. For example, they
remain true if $\Gamma$ is obtained from $\Gamma(R,h)$ by
replacing the two helical arcs of $\Gamma(R,h)$ with two smooth
curves in
$$
  \overline{H^+}\cap \{(x,y,z): x^2+y^2=R^2\}
$$
such that the angle function $\theta$  strictly increases from $0$
to $h$ on one curve and from $-h+\pi$ to $\pi$ on the other.  No
changes are required in the proofs.
\end{remark}

\stepcounter{theorem}\subsection{Vertical tangent planes along the
axis $Z$} \label{Zverticals}

Let $V$ be a vertical plane containing $Z$. Note that $M\cap V$
contains an interval in $Z$. If $M$ and $V$ are tangent at a point
$p\in Z$, then $M\cap V$ also contains a curve transverse to $Z$
that passes through $p$. Equivalently, $D\cap V$ contains a curve
with $p$ as one of its endpoints. Analyzing the curves of
intersection of $D$ and $V$ will allow us to understand the
distribution of points on $Z$ where $V$ and $M$ are tangent. This
analysis is carried out in Theorem~\ref{VtangentZ}.

We begin with a simple lemma.
\begin{lemma}\label{simplelemma}
The intersection of $\overline{D}$ with any plane in $\RR^3$
cannot contain a closed curve unless the curve passes through the
origin.
\end{lemma}

\begin{proof} Consider a closed curve $C$ in
         $\overline{D}$
not passing through the origin.  Since
  $\overline{D}\setminus\{0\}$
is an embedded disk, $C$ bounds a portion $D'$ of $D$. If a plane
contained $C$, then by the maximum principle that plane would
contain $D'$, and therefore by analytic continuation it would
contain all of $D$, which is impossible since $\Gamma$ is not a
planar curve.
\end{proof}

\begin{theorem}\label{VtangentZ}
Let $V$ be a plane that contains $Z$.  Let $W$ be a connected
component of $H^+\cap V$ that intersects $D$ transversely and that
does not contain any horizontal segment of $\Gamma(R,h)$ in its
boundary.  Then $W\cap D$ consists of at most three curves, and at
most four of the endpoints of those curves are on $Z$.
\end{theorem}

 \begfig
 %\hspace{0.3in}
\vspace{.2in} \centerline {
\includegraphics[width=3.0in]{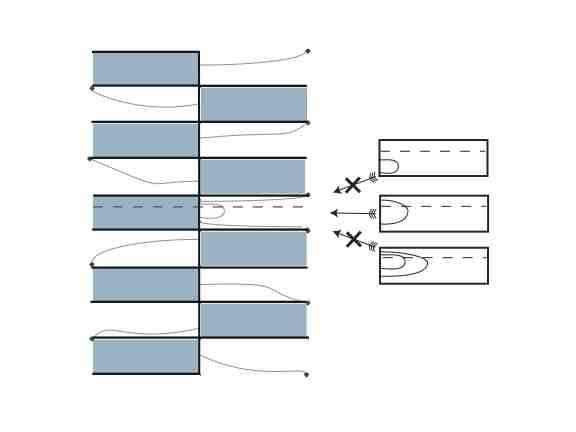}
} \vspace{-.5in}
 \hfill
 \begin{center}
 \vspace{0.3cm}
  \caption{}
   \label{figure:VerticalSlice}
   \vspace{.3cm}
    \parbox{4.0in}{
The possible intersections of a vertical plane $V$ and an embedded
minimal disk $D$, where
  $Z\subset V$,
  $D\subset H^+$,
and
  $\partial D = \Gamma(R,h)$.  The dotted horizontal
line is at the level $z=0$. The straight line segments represent
$H\cap V$.  The unshaded rectangles lie in $H^+\cap V$, the shaded
ones in $H^-\cap V$. The curves in $H^+\cap V$ indicate {\em
possible} intersections of $V$ and $D$. The component of
  $H^+\cap V$
that contains points at the level of $z=0$  is the only component
that can contain more than one intersection curve, and the only
component in which an intersection curve might conceivably have
both endpoints on $Z$.
  }
 \end{center}
 \endfig

\begin{proof} (See Figure~\ref{figure:VerticalSlice}.)
Note that $W=\sigma_\theta(S)$ for some $\theta\in \RR$, where $S$
is the halfstrip
 $$
    \{ (x,0,z): x>0,\, -\pi<z<0 \}
 $$
as in Section~\ref{s:angle function}.  Equivalently,
$$
  W = \{ p\in H^+: \theta(p) = \theta \}.
$$
By transversality, $D\cap W$ consists of smooth curves.

By Lemma~\ref{simplelemma}, none of the curves in $D\cap W$ is
closed. By elementary topology, the endpoints of such curves must
lie in the set
$$
  (\partial D)\cap \overline{W}  = I \cup P
$$
where $I$ is an interval in the $z$-axis, namely
$$
  \{ (0,0,z): z\in [-h,h]\cap [\theta-\pi,\theta]\}
$$
and where $P$ is the intersection of the helical portion of
$\Gamma$ with $\overline{W}$.

Note that $P$ consists of zero, one, or both (depending on
$\theta$) of the points
$$
  a(\theta) = (R\cos\theta, R\sin\theta, \theta)
$$
and
$$
   b(\theta) = (R\cos\theta, R\sin\theta, \theta-\pi).
$$
(See~\eqref{e:upper} and~\eqref{e:lower} in Section~\ref{s:angle
function}.)  Also note that each element of $P$ is indeed the
endpoint of exactly one curve of $D\cap W$. Thus at most two
curves in $D\cap W$ can have an endpoint in $P$.  These curves
together have at most two endpoints on the $z$-axis. Any other
curve $C$ of $D\cap W$ must have both endpoints on the $z$-axis.
Note one endpoint of $C$ must be on the positive $z$-axis and the
other on the negative $z$-axis, since otherwise $C$ together with
the segment joining its endpoints would violate
Lemma~\ref{simplelemma}. Furthermore, there cannot be a second
curve $C'$ from the positive $z$-axis to the negative $z$-axis,
since $C\cup C'$ together with the two segments joining their
endpoints (and not containing the origin) would then violate
Lemma~\ref{simplelemma}.
\end{proof}

If we assume that $\overline{W}$ is disjoint from the plane $z=0$,
then the method of proof of Theorem~\ref{VtangentZ} gives
considerably more:

\begin{proposition}\label{prop:moreVintersections}
 Suppose
$W=\sigma_\theta(S)=\{p\in H^+: \theta(p)=\theta\}$.
\begin{enumerate}
\item If $\overline{W}$ is disjoint from the plane $z=0$, then $D$
intersects $W$
  transversely.
\item If $\overline{W}$ is contained in the interior of the slab
  \[
     K =  \{ (x,y,z): |z| \le \max\{h, \pi\} \},
  \]
  and is disjoint from the plane $z=0$,
  then $D\cap W$ consists of a single smooth embedded curve with exactly one endpoint on $Z$.
\item If $\overline{W}$ is not contained in the slab $K$, then
$D\cap W$ is empty. \item If $\overline{W}$ is contained in the
slab $2\pi< z< h$ or in the slab $-h<z<-2\pi$,
  then $D\cap W$ is a single smooth embedded curve with one endpoint on $Z$, and that curve
  is a graph over a line segment in the plane $z=0$.
\end{enumerate}
\end{proposition}

\begin{proof} The hypothesis in assertion~(1) is equivalent to the condition that $\theta$ is not in the interval  $[0,\pi]$.  The transversality thus follows from Corollary~\ref{cor:vertical}.

In assertion~(2), the set $P$ in the proof of
Theorem~\ref{VtangentZ} consists of a singe point, and that proof
then shows that $D\cap W$ is a single curve with exactly one
endpoint on $Z$.

In assertion~(3), the set $P$ is empty, from which it follows that
$D\cap W$ is also empty.

In assertion~(4), $D\cap W$ consists of a single smooth curve by
assertion~(2). The hypothesis implies that $\theta$ is not
contained in the interval $[-\pi,2\pi]$.  By Theorem~\ref{vertical
thm 1}, that curve has no vertical tangents and thus is a graph
over its projection to the plane $z=0$.
\end{proof}

\stepcounter{theorem}\subsection{An upper bound on the number of
tangent planes parallel to a vertical plane $V$}\label{upperbound}

Let $D$ be a $\rho_Y$-invariant minimal embedded disk in $H^+$
with boundary $\Gamma=\Gamma(R,h)$, let
 $$
   M = \overline{D}\cup \rho_Z(\overline{D}),
 $$
and let $N$ be the corresponding $\sigma_{2h}$-invariant surface.
We combine the results of the two previous sections to get a local
bound on the number of points of $N$ at which that the tangent
plane is parallel to a given vertical plane.  The estimate is
independent of $R$ and $h$.

\begin{theorem}\label{Vestimate}
Let $V$ be a plane containing $Z$, and let $\Pi$ be an open
horizontal slab of thickness $\pi$:
$$
    \Pi =\{(x,y,z): a < z < a+\pi \}.
$$
Then $\Pi$ contains at most sixteen points of $N$ at which the
tangent plane is parallel to $V$.
\end{theorem}

The number sixteen is certainly not optimal, but for our purposes
any finite number would suffice.

\begin{proof}
Consider first a plane $V$ containing $Z$ that is generic in the
following sense: $V$ intersects $N\setminus Z$ transversely and
$V$ does not contain any of the countably many horizontal line
segments in $N$. Now $N$ is made up of congruent copies of $M$,
which is in turn made up of two copies of $\overline{D}$. Note
that $\Pi$ contains points from at most two copies of $M$, and
therefore from at most four copies of $\overline{D}$. Now
$\Tan_pN$ is not vertical at any boundary point of $N$ by the
boundary maximum principle. Also, $V$ does not contain any of the
horizontal edges of the copies of $\overline{D}$.  Thus if
$\Tan_pN$ is parallel to $V$, then $p$ lies either on the $z$-axis
or else in the interior of one of the four copies of
$\overline{D}$. By Theorem~\ref{VtangentZ} (see also the
discussion in~\ref{Zverticals}), there are at four such points on
$Z$ in each copy of $M$.  By Corollary~\ref{cor:vertical}, there
are at most two such points not on $Z$ in each copy of $D$.   Thus
there are at most $2\times 4 + 4\times 2$ or sixteen such points
in $\Pi$.

By openness of the Gauss map, the number of such points is a lower
semicontinuous function of $V$.  Thus the bound for arbitrary $V$
follows from the bound for generic $V$.
\end{proof}

\section{Uniform estimates}\label{section:uniform_est}
\stepcounter{theorem}\subsection{A uniform curvature estimate}

Consider a $\rho_Y$-invariant minimal embedded disk $D$ in $H^+$
with boundary $\partial D=\Gamma =\Gamma(R,h)$, the curve
specified in Section~\ref{s:defs}.  Extend $D$ by Schwarz
reflection in $Z$ to get a minimal embedded genus-one surface $M$
with $\partial M =\partial H_{R,h}$. As observed in
Section~\ref{helicoidguide}, this boundary consists of a top line
segment, a bottom line segment, and two helical arcs. By repeated
Schwarz reflection in the top and bottom line segments, we get a
smooth embedded minimal surface $N$ invariant under the screw
motion $\sigma_{2h}$.  The boundary of $N$ consists of two
helices.

In Section~\ref{section:main thms}, we will obtain complete
nonperiodic (or periodic) genus-one helicoids by letting $R$ and
$h$ (or just $R$) tend to infinity. In this section, we prove the
estimates that allow us to control passage to the limit.

\begin{lemma}\label{simple connectivity lemma}
Suppose $D$ is a $\rho_Y$-invariant minimal embedded disk in $H^+$
with boundary $\partial D=\Gamma=\Gamma(R,h)$.  Let
 $$
     M = \overline{D}\cup \rho_Y(\overline{D})
 $$
and let $N$ be the $\sigma_{2h}$-invariant surface obtained from
$M$.
\begin{enumerate}
 \item If $C$ is a closed curve in $N$ that does not intersect any
 straight line or straight line segment contained in $N$, then $C$
 is contractible in $N$.
 \item Any pair of disjoint homologically nontrivial embedded closed
 curves in $M$ must bound an annulus in $M$.
 \item There is a unique shortest homotopically nontrivial curve
 $\alpha$  in $\overline{D}$.
 It is a smooth closed geodesic in $N$ that is bisected by $0$
 together with the unique point of $Y\cap D$.
\end{enumerate}
\end{lemma}

\begin{proof}
Note that if we remove the straight line segments from $N$, we are
left with a disjoint union of pieces, each of which is congruent
to $D$ and therefore is simply connected.  In assertion (1), $C$
lies entirely in one of those pieces and is therefore contractible
in it.

In Section~\ref{subsection:replacing}, we showed that $M$ is
topologically a once-punctured torus. Assertion (2) follows by
standard, elementary topology.

Let $\alpha$ be a shortest curve that is a generator for
$\pi_1(\overline{D},0)\cong \ZZ$. (The fundamental group is
infinite cyclic since $\overline{D}$ is topologically a closed
disk with two boundary points identifed.) Note that $N$ has
geodesically convex boundary because the curvature vector at each
point of its bounding helices points toward the $z$-axis.  Since
$\overline{D}$ is bounded by portions of $\partial N$ together
with geodesics in $N$, the curve $\alpha$ does not touch $\partial
D$ except at $0$. Thus $\alpha$ is a smoothly embedded geodesic in
$\overline{D}$ that starts and ends at $0$.  Since $D$ has
nonpositive curvature, $\alpha$ is a unique.   Therefore it is
$\rho_Y$ invariant, so its midpoint $p$ is a fixed point of
$\rho_Y$, namely the unique point of $Y\cap D$. (See
Proposition~\ref{vertical thm 1}.) The $\rho_Y$ invariance also
implies that the two components of $\alpha\setminus\{0,p\}$ are
related by $\rho_Y$. This implies that $\alpha$ does not have a
corner at $0$, but rather forms a smooth closed geodesic in $N$.
\end{proof}

\begin{remark}\label{remark:simple connectivity lemma}
Note that conclusions~(1) and~(2) of Lemma~\ref{simple
connectivity lemma} are preserved under smooth convergence.  That
is, if $M_i$ and $N_i$ are smooth minimal surfaces satisfying the
first two conclusions of the lemma, and if the $M_i$ and $N_i$
converge smoothly to limits $M$ and $N$, then $M$ and $N$ also
satisfy those two conclusions.
\end{remark}

The proofs of our next two results rely strongly on the following
theorem of Mo and Osserman \cite{mos1}, which extends earlier work
of Osserman \cite{oss}, Xavier \cite{xa1}, Sa~Earp-Rosenberg
\cite{er3}, and Fujimoto \cite{fu1}.  (The Sa~Earp-Rosenberg
result is also strong enough for our purposes.)

\begin{theorem}[Mo-Osserman]\label{Mo-Osserman}
    If $\Nn$ is a complete minimal surface in $\RR^3$ with Gauss map
$g:\Nn\to {\mathbf S}^2$ and if the set
$$
 \{ v\in {\mathbf S}^2: \text{$g^{-1}(v)$ is finite} \}
$$
contains five or more points, then $\Nn$ has finite total
curvature.
\end{theorem}

We now give our main curvature estimate:

\begin{theorem}\label{thm:MainCurvatureEstimate}
 There are finite constants $R_0$ and $K$ with
the following property.  Suppose $D$ is an $\rho_Y$-invariant
minimal embedded disk in $H^+$ with boundary $\Gamma(R,h)$, where
$R\ge R_0$ and $h\ge \pi/2$.  Then
$$
 B(D,p)\le K.
$$
\end{theorem}
Here $B(D,p)$ is the norm of the second fundamental form of $D$ at
$p$.

The hypothesis $h\ge \pi/2$ could be removed, since one can show
that there is no such disk $D$ if $h<\pi/2$.  However, we only
require the theorem for $h\ge \pi/2$.

\begin{proof}
Suppose the theorem is false. Then there is a sequence of examples
$D_i$ with $\partial D_i=\Gamma(R_i,h_i)$ and a sequence of points
$p_i\in \overline{D}_i$ such that
  $B(D_i,p_i)$ and $R_i$ tend to infinity and such that $h_i\ge \pi/2$.

Let
 $$
   M_i = \overline{D}_i\cup \rho_Z(\overline{D}_i),
 $$
let $N_i$ be the screw-motion-invariant surface obtained from
$M_i$, and let $C_i$ be the solid cylinder of radius $R_i$ about
the $z$-axis.

We may suppose that $p_i$ has been chosen in $\overline{D}_i$ to
maximize $B(\overline{D}_i,p_i)=B(N_i,p_i)$.  (The maximum exists
because $N_i$ is smooth and $\overline{D}_i$ is a compact subset
of $N_i$.)  It follows that
$$
  \max_{p\in N_i}B(N_i,p) = B(N_i,p_i).
$$

Translate $N_i$, $M_i$, and $C_i$ by $-p_i$, and then dilate by
$B(N_i,p_i)$ to get $N_i'$, $M_i'$, and $C_i'$. By passing to
subsequences we may assume that $N_i'$, $M_i'$, and $C_i'$
converge as sets to limits $N'$, $M'$, and $C'$.

Note that
\begin{equation}\label{e:bi'=1}
  \max_{p\in N_i'} B(N_i',p) = B(N_i',0) = 1.
\end{equation}

Also, either $\partial N_i'$ converges to the empty set or else it
converges smoothly to a horizontal line.  (There is just one line
since $h_i\ge \pi/2$ and since the dilation factor $B(D_i,p_i)$
tends to infinity. The line is horizontal since $R_i\to\infty$.)
This together with~\eqref{e:bi'=1} implies that the convergence
$N_i'\to N'$ is smooth.  Thus the limit is a smooth, embedded
minimal surface and $\partial N'$ is either the empty set or a
horizontal line.

Let
\begin{equation*}
   \Nn = \begin{cases} N' &\text{if $\partial N'=\emptyset$} \\
                       N'\cup \rho_{\partial N'}(N')
                          &\text{if $\partial N'\ne \emptyset$}.
         \end{cases}
\end{equation*}
If $\partial N'$ is a line, then $N'$ lies in a halfspace (namely
$C'$) whose boundary plane contains that line, and thus $\Nn$ is
embedded. Of course if $\partial N'$ is empty, then $\Nn$ is also
embedded.  Either way, $\Nn$ is a complete embedded minimal
surface.

By Theorem~\ref{Vestimate}, $N'$ has the following property: if
$V$ is a vertical plane, then there are at most a finite number of
points of $N'$ at which the tangent plane is parallel to $V$. (If
this is not clear, recall that the $h_i$ are bounded below and
that the dilation factor $B(D_i,p_i)$ tends to infinity.)  Note
that if $L$ is a horizontal line, then $\rho_L(N')$ has the same
property. Consequently $\Nn$ also has this property. In other
words, the set
$$
   \{ v\in {\mathbf S}^2: \text{$g^{-1}(v)$ is finite} \}
$$
contains a great circle, where $g:\Nn\to {\mathbf S}^2$ is the
Gauss map.  Hence by the Mo-Osserman theorem (\ref{Mo-Osserman}),
$\Nn$ has finite total curvature.

We now know that $\Nn$ is a complete embedded minimal surface of
finite total curvature. By~\eqref{e:bi'=1},
 \begin{equation}\label{e:b=1}
    B(\Nn,p')=1,
 \end{equation}
so $\Nn$ is not flat.  Therefore $\Nn$ has a catenoidal end.
Intersecting the end with a suitable plane parallel to the end, we
see that $\Nn$ contains a planar closed curve that does not
intersect any straight line segment contained in $\Nn$.  By
assertion~(1) of Lemma~\ref{simple connectivity lemma}, that curve
bounds a disk in $\Nn$. By the maximum principle, the disk must
lie in the plane containing its boundary. But then by analyticity,
all of $\Nn$ must be planar, contradicting~\eqref{e:b=1}.
\end{proof}

The following theorem will let us conclude that the complete
surfaces we construct are nearly horizontal away from the
$z$-axis.

\begin{theorem}\label{translates}
Suppose $D_i$ is a sequence of embedded, $\rho_Y$-invariant
minimal disks in $H^+$ with $\partial D_i=\Gamma(R_i,h_i)$, where
$R_i\to\infty$ and
 $
     h_i\to h\in (\pi/2,\infty].
 $
Let $D_i'$ be the result of translating $D_i$ by $-p_i$, where
$p_i\in D_i$ is a sequence of points such that
\begin{equation}\label{e:translates thm hypothesis}
  \dist(p_i,Z\cup \partial C_i)\to \infty.
\end{equation}
Here $\dist(p_i,\cdot)$ denotes intrinsic distance in $D_i$ and
$C_i$ denotes the solid cylinder of radius $R_i$ about $Z$.

After passing to a subsequence, the $D_i'$ converge smoothly to a
limit $D'$.  Let $\Sigma$ be the component of $D'$ containing $0$.
Then
\begin{enumerate}
 \item $\Sigma$ is a horizontal plane, or
 \item $\Sigma$ is a horizontal halfplane, or
 \item $\partial \Sigma$ consists of two lines parallel to the $x$-axis.
\end{enumerate}
The third case can occur only if $h=\pi$, the intrinsic
$\dist(p_i,X)$ from $p_i$ to $X$ is bounded, and the length of the
shortest closed geodesic in $\overline{D_i}$ tends to infinity.
\end{theorem}

\begin{remark} One can prove (using a slight generalization of
Proposition~\ref{halfslab} below) that in case (3), the surface
$\Sigma$ must be a flat strip.  However, we will not need that
fact. Our main construction (Section~\ref{section:main thms}) uses
disks with the annular intersection property, and
Theorem~\ref{thm:BoundOnGeodesic} below implies that case~(3) does
not occur for disks $D_i$ having that property.
\end{remark}

\begin{proof} By the convex hull property, $D_i$ lies in
the solid cylinder $C_i$.  By
Theorem~\ref{thm:MainCurvatureEstimate}, we may assume (after
passing to a subsequence) that the $D_i'$ converge smoothly to
$D'$ and that the corresponding translates of the $C_i$ and of $H$
converge to a limits $C'$ and $H'$, respectively.

Note that $H'$ is either a helicoid or else a union of horizontal
planes according to whether the Euclidean distance from $p_i$ to
$Z$ stays bounded or tends to infinity.  Similarly, since
$R_i\to\infty$, the limit $C'$ is either a halfspace bounded by a
vertical plane or else all of $\RR^3$ according to whether the
Euclidean distance from $p_i$ to $\partial C_i$ stays bounded or
tends to infinity.

By Corollary~\ref{cor:vertical}, for each vertical plane $V$,
there are at most two points of $\Sigma$ at which the tangent
plane is parallel to $V$.

Consider first the case $\partial \Sigma=\emptyset$.  By the
Mo-Osserman Theorem~\ref{Mo-Osserman}, $\Sigma$ has finite total
curvature. Also, it is complete, embedded, and simply connected.
Therefore it is a plane.

Since $H$ and $D_i$ are disjoint, $H'$ and $\Sigma$\, cannot cross
each other (i.e., contain points of transverse intersection). Now
$H'$ is either a helicoid or a union of horizontal planes. Since
the plane $\Sigma$ crosses every helicoid, $H'$ must be a union of
parallel planes.  Since $\Sigma$ does not cross $H'$, the plane
$\Sigma$ must also be horizontal.

Note that the horizontal plane $\Sigma$ is contained in $C'$, so
$C'$ must be all of $\RR^3$ (rather than a halfspace bounded by a
vertical plane.)  Since $H'$ is a union of horizontal planes and
$C'$ is all of $\RR^3$, it follows, as explained above, that the
Euclidean distance from $p_i$ to $Z\cup \partial C_i$ tends to
infinity in the case $\partial \Sigma=\emptyset$.

Now suppose $\partial \Sigma$ is not empty. By~\eqref{e:translates
thm hypothesis}, the boundary $\partial \Sigma$ consists of one or
more horizontal lines corresponding to the horizontal radial
segments in $\partial D_i$.  Since intrinsic distance in
$\overline{D_i}$ and Euclidean distance coincide on line segments
in $\overline{D}_i\setminus\{0\}$, the
hypothesis~\eqref{e:translates thm hypothesis} implies that the
Euclidean distance from $p_i$ to $Z\cup\partial C_i$ tends to
infinity. This in turn implies that $H'$ is a union of horizontal
planes.

Suppose that $\partial \Sigma$ consists of a single horizontal
line $L$.  Let
$$
  \Sigma^*=\Sigma\cup L \cup \rho_L\Sigma.
$$
Then (just as before) $\Sigma^*$ is a complete, embedded, simply
connected minimal surface of finite total curvature.  Thus
$\Sigma^*$ is a plane, so $\Sigma$ is a half-plane.  Since $D_i$
lies in $H^+$, $\Sigma$ must lie in the closed region between two
successive planes in $H^+$.  Thus $\Sigma$ is horizontal.

Finally, suppose $\partial \Sigma$ consists of more than one
horizontal line.  Since $D_i$ lies in $H^+$, $\Sigma$ lies in the
region between two successive planes in $H^+$.  That is, $\Sigma$
lies in a horizontal slab $\Omega$ of thickness $\pi$.  It follows
that $\partial \Sigma$ consists of exactly two horizontal lines,
one in each component of $\partial \Omega$.  The vertical distance
between those two lines is $\pi$.

On the other hand, the horizontal segments in $\partial
D_i=\Gamma(R_i,h_i)$ lie in the planes $z=-h_i$, $z=0$, and
$z=h_i$.  Thus the vertical distance between two lines in
$\partial \Sigma$ is either $h$ or $2h$. Thus $h=\pi$ or $2h=\pi$.
Since $h>\pi/2$, this means $h=\pi$.

Since $h_i\to \pi$, the horizontal edges of $\partial D_i$
converge to the positive and negative portions of the $x$-axis and
to the rays
 $
  \{ (x,0,\pi): x\le 0\}
 $ and
 $
 \{ (x,0,-\pi): x \ge 0\}.
 $
Thus the two lines of $\partial \Sigma$ are parallel, and (after
passing to a subsequence if necessary) they are limits either:
\begin{enumerate}
 \item of $X^-$ and of the top horizontal edge $T_i$, both translated by $-p_i$, or
 \item of $X^+$ and the bottom horizontal edge $B_i$, both translated by $-p_i$.
\end{enumerate}
Of course in either case $\dist(p_i,X)$ must be bounded.  Note in
case (1), $\dist(p_i,T_i)$ is also bounded, and in case (2)
$\dist(p_i,B_i)$ is bounded.

Without loss of generality, assume we are in case (1), so that
$\dist(p_i,X^-)$ and $\dist(p_i,T_i)$ are bounded.

Let $\alpha_i$ be the shortest closed curve in $\overline{D_i}$.
Then $\alpha_i$ divides $\Gamma(R_i,h_i)$ into two components. One
component consists of the top edge $T_i$, the edge in postive
$x$-axis, and a helical arc joining them.  The other component
consists of the bottom edge $B_i$, the edge in the negative
$x$-axis, and a helical arc joining them.  In particular, the top
horizontal edge and the edge in the negative $X$-axis belong to
different components of $\Gamma(R_i,h_i)\setminus \alpha_i$ and
thus in different components of $D_i\setminus \alpha_i$.  Hence
either the shortest curve from $p_i$ to $T_i$ or the shortest
curve from $p_i$ to $X^-$ must cross $\alpha_i$. Therefore the
union of those two curves with $\alpha_i$ contains a path joining
$p_i$ to the origin. Thus
$$
  \dist(p_i,0) \le \Length(\alpha_i) + \dist(p_i,T_i) + \dist(p_i,X^-).
$$
Since the left hand side tends to infinity and since the second
and third terms on the right are bounded, $\Length(\alpha_i)$ must
tend to infinity.
\end{proof}

In proving the theorem, we also proved

\begin{corollary}\label{cor:properness}
If the intrinsic distance from $p_i$ to
 $Z\cup \partial C_i$
tends to infinity, then the Euclidean distance from $p_i$ to
$Z\cup \partial C_i$ also tends to infinity.
\end{corollary}

The following special case of Theorem~\ref{translates} will be
used in~\ref{thm:BoundOnGeodesic}:

\begin{corollary}\label{cor:translates}
If $\dist(p_i,X\cup Z\cup \partial C_i)\to \infty$, then
$\Tan_{p_i}D_i$ converges to a horizontal plane.
\end{corollary}

\begin{proof} Since $\dist(p_i, Z\cup\partial C_i)\to\infty$,
assertion (1), (2), or (3) of Theorem~\ref{translates} must hold.
Since $\dist(p_i,X)\to\infty$, assertion (3) does not hold.
\end{proof}

\stepcounter{theorem}\subsection{A uniform bound on the length of
the closed geodesic in $\overline{D}$}

Consider a $\rho_Y$-invariant, minimal embedded disk $D$ in $H^+$
with boundary $\Gamma(R,h)$ (for some $R$ and $h$), where
$\Gamma(R,h)$ is the curve defined in Section~\ref{s:defs}.  The
next theorem establishes a uniform estimate for the length of the
shortest closed geodesic in $\overline{D}$, provided the disk $D$
has the annular intersection property~\ref{def:annular
intersection property}. (This is the only estimate in the paper
that depends on the annular intersection property.) We will use
this result to show that the genus-one surfaces
$M=\overline{D}\cup \rho_Z(\overline{D})$ have genus-one limits as
$R\rightarrow \infty$.

\begin{theorem}\label{thm:BoundOnGeodesic}
Suppose $D_i$ is a sequence of  $\rho_Y$-invariant minimal
embedded disks in $H^+$ with boundary $\Gamma(R_i,h_i)$ and with
the annular intersection property.  Suppose also that
$$
  h_i \ge \eta > \pi/2
$$
and that $R_i\to \infty$. Then the length of the shortest closed
geodesic in $\overline{D_i}$ is bounded above.
\end{theorem}

\begin{proof} Let $\alpha_i$ be the shortest closed geodesic in $\overline{D}_i$.
By assertion~(3) of Lemma~\ref{simple connectivity lemma},
$\alpha_i$ contains the unique point $p_i$ in $D_i\cap Y$, and the
length of $\alpha_i$ is twice the intrinsic distance from $p_i$ to
$0$. Thus it suffices to bound $\dist(p_i,0)$ above.

Let $q_i$ be the point in $X\cup Z$ that is closest in intrinsic
distance to $p_i$. Then
 \begin{align*}
  \dist(p_i,0) &\le \dist(p_i,q_i) + \dist(q_i,0)
                     \\
               &= \dist(p_i,q_i) + |q_i|
                     \\
               &\le \dist(p_i,q) + |p_i - q_i|
                     \\
               &\le 2 \dist(p_i,q_i)
 \end{align*}
(Here $|q_i|\le |p_i-q_i|$ because $p_i\in Y$ and $q_i\in X\cup Z$
are orthogonal.) Thus
 \begin{equation}\label{e:dist to XunionZ}
    \dist(p_i,0) \le 2\dist(p_i,X\cup Z).
 \end{equation}
The tangent plane to $D_i$ at $p_i$ is vertical by
Proposition~\ref{vertical thm 1}.
 Thus by
Corollary~\ref{cor:translates}, the sequence
\begin{equation*}
            \dist(p_i, X\cup Z \cup \partial C_i)
\end{equation*}
is bounded, which implies by~\eqref{e:dist to XunionZ} that the
sequence
$$
   \dist(p_i, \{0\}\cup\partial C_i)
$$
is bounded.  Hence it suffices to prove that
\begin{equation}\label{e:far from cylinder}
   \dist(p_i, \partial C_i) \to \infty.
\end{equation}

Suppose, on the contrary, that $\dist(p_i,\partial C_i)$ is
bounded.  Translate $D_i$, $H$, and $C_i$ by $-p_i$ to get $D_i'$,
$H_i'$, and $C_i'$.  Since $\dist(p_i,\partial C_i)$ is bounded
and $R_i\to\infty$, the $C_i'$ converge (after passing to a
subsequence) to a halfspace $C'$ of the form
 $$
   C' = \{(x,y,z): y\le a\}.
 $$
It follows that the $H_i'$ converge smoothly to a limit $H'$
consisting of the horizontal planes on which $z$ is an odd
multiple of $\pi/2$.  By the curvature estimate in
Theorem~\ref{thm:MainCurvatureEstimate}, we may assume that the
$D_i'$ converge smoothly to a limit surface. Let $D'$ be the
connected component of that limit surface containing the origin.
Then $D'$ lies in one of the components of $C'\setminus H'$,
namely the halfslab
\begin{equation}\label{e:halfslab}
  \{ (x,y,z): \text{$y\le a$ and $|z|<\pi/2$} \}.
\end{equation}
Note that $\partial D'$ consists of the two straight line edges of
the halfslab.  (By the $\rho_Y$ symmetry of $D'$, the boundary
$\partial D'$ cannot be just one of the two lines.)

By Corollary~\ref{cor:properness}, $\overline{D'}$ is properly
embedded.

We have shown that $\overline{D'}$ is a simply connected properly
embedded minimal surface in the half-slab~\eqref{e:halfslab} and
that $\partial D'$ consists of the two edges
of~\eqref{e:halfslab}.  The only such minimal surface is the
vertical strip $S$ bounded by those two edges (see Proposition
\ref{halfslab} below), so $D'=S$. However, we will show that
$D'=S$ contradicts the annular intersection property of the $D_i$.

Let $A'$ be a catenoid that intersects the planes $z=\pm\,\pi/2$
in a pair of circles in the region $y<a$.  Since $D_i$ intersects
each annulus in $\Aa(\Gamma_i)$, it follows that $D'$ intersects
$A'$. But the strip $D'=S$ does not intersect $A'$, so we have a
contradiction.

(In case it is not clear, we spell out in more detail why $D'$
must intersect $A'$.  Let $A_i'$ be the component of $A'\setminus
H_i'$ that crosses the $xy$-plane.  Translate $A_i'$ by $p_i$ to
get $A_i$.  Because of the smooth convergence, when $i$ is
sufficiently large, $A_i$ will be in $\Aa(\Gamma_i)$, and so $D_i$
will intersect $A_i$. Thus $D_i'$ intersects $A_i'$ and hence,
passing to the limit, $D'$ intersects $A'$.)
\end{proof}

\begin{remark}\label{remark:aip not needed}
Note that the annular intersection property was only used to
prove~\eqref{e:far from cylinder}.  Thus, even without assuming
the annular intersection property, the theorem applies to any
sequence $D_i$ for which \eqref{e:far from cylinder} holds.
\end{remark}

\begin{proposition}\label{halfslab}
Suppose $U\subset\RR^3$ is a half-slab bounded by two horizontal
half-planes and an infinite strip $S$.  Suppose $M$ is a simply
connected, properly embedded minimal surface in $\overline{U}$
with $\partial M=\partial S$.  Then $M=S$.
\end{proposition}

\begin{proof}
We may suppose that $
  U = \{(x,y,z): \text{$y< 0$ and $|z|<h$} \}
$ and thus that the strip is $
  S = \{(x,y,z): \text{$y=0$ and $|z|<h$}\}.
$

Since $M$ is connected, it contains an embedded path $\gamma$
joining $(0,0,-h)$ to $(0,0,h)$.  Since $M$ is simply connected,
$\gamma$ divides $M$ into two components.  One of those
components, which we will denote $M^+$, has boundary consisting of
$\gamma$ together with $\{(x,0,\pm h)\in \partial S: x\ge 0\}$.

(If this is not clear, note that if $M\ne S$, then $M\cup S$
bounds a region in $\RR^3$ and so is orientable. Thus $M\cup S$ is
topologically an annulus (rather than a M\"obius strip), and thus
$\gamma$ together with the segment joining its endpoints divides
$M\cup S$ into two components.)

Let $\BB$ be a ball centered at the origin and containing
$\gamma$. Let $\Cc$ be the set of catenoids $C$ with vertical axes
of symmetry such that $C$ is disjoint from $
   \BB\cup \{(x,y,z)\in S: x \ge 0\}.
$ Note that $\Cc$ is a connected family, and that there are
catenoids in $\Cc$ that are disjoint from $M+$.  (Take any $C\in
\Cc$ disjoint from $\overline{U}$.)  Thus by the maximum
principle, all the catenoids in $\Cc$ are disjoint from $M^+$.
These catenoidal barriers force $M^+$ to be arbitrarily close to
$S$ near $\infty$:
 $$
   \lim_{(x,y,z)\in M^+,\, x\to\infty} y = 0.
 $$
Likewise for $(x,y,z)\in M^-=M\setminus M^+$, we see that $y\to 0$
as $x\to -\infty$.  Thus $|y|$ achieves its maximum on $M$ at a
finite point $p$.  By the maximum principle, this maximum value
must be $0$.
\end{proof}

\section{The Main Theorems}\label{section:main thms}

\begin{theorem}\label{thm:main theorem 1}
There exists a complete, properly embedded minimal surface $M$ in
$\RR^3$ such that
\begin{enumerate}
 \item\label{mainthm1:McapH} $M\cap H=X\cup Z$.
 \item\label{mainthm1:SimplyConnectedPieces}
  Each of the two components of $M\setminus H$ is simply
 connected.
 \item\label{mainthm1:torus}
            $M$ is topologically a once-punctured torus.
 \item\label{mainthm1:horizontal} $\Tan_pM$ converges to a horizontal plane as
 $\dist(p,Z)\to\infty$.
 \item\label{mainthm1:VerticalTangents}
  The points of $M\setminus Z$ with vertical tangent planes
 lie in a cylinder $\BB(0,R)\times [-2\pi, 2\pi]$.
 \item\label{mainthm1:conformaltype} $M$ is conformally a once-punctured torus.
 \item\label{mainthm1:level0}  The level set $M\cap \{z=0\}$  consists of  $X$ together with a
 smooth, simple closed curve that intersects $X$ in exactly two points.
 \item\label{mainthm1:levelc}  For each $c\ne 0$, the level set $M\cap \{z=c\}$
consists of a single smooth, nonclosed curve.
\item\label{mainthm1:asymptotic} $M$ is asymptotic to $H$ at
infinity.
\end{enumerate}
\end{theorem}

\begin{proof}
Choose sequences $R_i$ and $h_i$ tending to infinity.   By
Propositions~\ref{prop:aip} and~\ref{existence3}, for all
sufficiently large $i$, the curve $\Gamma(R_i,h_i)$ bounds a
$\rho_Y$-invariant minimal embedded disk $D_i$ in $H^+$ with the
annular intersection property. Let
$$
  M_i = \overline{D}_i\cup \rho_Z(\overline{D}_i).
$$
By Theorem~\ref{thm:MainCurvatureEstimate}, the curvatures of the
$M_i$ are uniformly bounded, so (by passing to a subsequence) we
may assume that the $M_i$ converge to a complete, embedded minimal
surface $M$.

Now
$$
  \partial (M_i\cap H^+) = \partial D_i = \Gamma(R_i,h_i) \to X\cup Z.
$$
It follows that
$$
  \partial (M\cap H^+) = X\cup Z
$$
which implies assertion~\eqref{mainthm1:McapH}.

Since $M\cap H^+$ is the limit of the simply connected minimal
surfaces $M_i\cap H^+=D_i$, it must also be simply connected.
Similarly $M\cap H^-$ must be simply connected. This proves
assertion~\eqref{mainthm1:SimplyConnectedPieces}.

Furthermore, $M$ is proper by Corollary~\ref{cor:properness}.

By assertion~(3) of Lemma~\ref{simple connectivity lemma}, $M_i$
contains a simple closed geodesic $\alpha_i$ such that $0\in
\alpha_i$ and such that $\alpha_i\setminus\{0\}$ lies in $D_i$.
Thus $\alpha_i'=\rho_Z(\alpha_i)$ is another closed geodesic, and
$\alpha_i$ and $\alpha_i'$ intersect transversely at the origin
and nowhere else. (They intersect only at the origin because
$\alpha_i\setminus \{0\}$ lies in $D_i\subset H^+$ and
$\alpha_i'\setminus \{0\}$ lies in $\rho_Z(D_i)\subset H^-$.)

The lengths of the geodesics $\alpha_i$ and $\alpha_i'$ are
bounded above by Theorem~\ref{thm:BoundOnGeodesic}, and they are
bounded below since the curvatures of the $M_i$ are uniformly
bounded (Theorem~\ref{thm:MainCurvatureEstimate}). Thus (after
passing to a subsequence if necessary) the $\alpha_i$ and
$\alpha_i'$ converge to closed geodesics $\alpha$ and $\alpha'$ in
$M$ that intersect transversely at the origin. Thus $M$ has genus
at least one. By part (2) of Lemma~\ref{simple connectivity lemma}
(see also Remark~\ref{remark:simple connectivity lemma}), $M$ has
genus at most one. Thus $M$ has genus exactly one.  Assertion~(2)
of Lemma~\ref{simple connectivity lemma} also implies that $M$ has
exactly one end. Thus $M$ is topologically a once-punctured torus.

Assertion~\eqref{mainthm1:horizontal} follows from
Theorem~\ref{translates}.

Assertion~\eqref{mainthm1:VerticalTangents} follows from
assertion~\eqref{mainthm1:horizontal} together with
Corollary~\ref{cor:vertical}.

To prove the remaining assertions, it is convenient first to prove
the following:
\begin{claim} Let $M^+=M\cap H^+$.
 \begin{enumerate}[\upshape (i)]
 \item\label{claim:OneComponent}
 If $c\ge 0$, then
 $
    M^+ \cap \{z>c\}
 $
 has exactly one connected component.
 \item\label{claim:NoHorizontalPlanes}
      If $\Tan_pM$ is horizontal, then $p\in X\setminus\{0\}$.
 \item\label{claim:SingleCurve}
 $M^+\cap \{z=0\}$ consists of a single smooth embedded curve.
  \item\label{claim:Dichotomy}
 Either $M\cap \{z=0\}$ consists of three connected components each of which is a smooth
 embedded curve, or it consists of $X$ together with
a smooth, simple closed curve that crosses $X$ exactly twice.
\end{enumerate}
\end{claim}

\begin{proof}[Proof of claim]
To prove~\eqref{claim:OneComponent}, let $C$ be the component of
$M^+\cap\{z>c\}$  that contains
\[
  \{ (0,0,z): z > c\}
\]
in its boundary, and let
\[
  C' = (M^+\cap \{z>c\}) \setminus C.
\]
 If $p\in M^+\cap \{z>c\}$ and if $\theta(p)>c+2\pi$ (where $\theta$ is the angle function defined
 in Section~\ref{s:angle function}), then by
assertion~(4) of Proposition~\ref{prop:moreVintersections}, there
is a curve in $M^+\cap\{z>c\}$ that contains $p$ and that has an
endpoint on $Z$.  Thus such a point $p$ must lie in the component
$C$.  This shows that the function $\theta$   is bounded above on
$C'$.
  Thus the function $z$ is also bounded above on $C'$. (Recall that $\theta(x,y,z)$ and $z$ differ by at most $\pi$.)

Now $C'$ is a minimal surface that is properly embedded in
$\{z>c\}$, $C'$ is contained in a slab, and $\partial C^{'}\subset
\{z=c\}$. A version of the halfspace theorem
\cite{HoffmanMeeksHalfspace} states that if a connected, properly
immersed minimal surface $\Sigma$    lies in a slab, and $\partial
\Sigma$ (if nonempty) lies on one boundary face of the slab, then
$\Sigma$ is a subset of a plane. Hence $C'$ is a union of
horizontal planes. But $C'$ is contained in $H^+$, so $C'$ must be
empty. This proves~\eqref{claim:OneComponent}.

To prove~\eqref{claim:NoHorizontalPlanes}, suppose the horizontal
plane $\{z=c\}$ is tangent to $M$ at some point $p$. By the
$\rho_Y$ symmetry, we may assume without loss of generality that
$c\ge 0$. If $p$ were in $M^+$, then by a theorem of Rado
(\cite{RadoBook}, III.7 or \cite{os1}, Lemma~7.5), the tangent
plane $\{z=c\}$ would divide $M^+$ into four of more components,
at least two of which would  lie in the region $z>c$,
contradicting~\eqref{claim:OneComponent}.  Thus $p$ does not lie
in $M^+$. By the same argument (or by $\rho_Z$ symmetry), it also
cannot lie in $M^-$, the other component of
  $M\setminus H$.
Thus by assertion~\eqref{mainthm1:McapH}, $p\in X\cup Z$.  Since
$Z\subset M$, the tangent plane to $M$ at every point of $Z$ is
vertical.  Thus $p\in X\setminus\{0\}$. This
proves~\eqref{claim:NoHorizontalPlanes}.

To prove~\eqref{claim:SingleCurve}, note that if $M^+\cap\{z=0\}$
 contained more than a single embedded curve, then
$M^+\setminus \{z=0\}$ would have more than two components.  By
the $\rho_Y$ symmetry, it would have more than one component in
the halfspace $\{z>0\}$.  But to according
to~\eqref{claim:OneComponent} (with $c=0$), there is only one such
component.

To prove~\eqref{claim:Dichotomy}, let $S$ be the curve whose
existence is given by statement~\eqref{claim:SingleCurve} of the
claim.   If $\overline{S}$ has no endpoints, then it follows from
statement~\eqref{claim:SingleCurve}  that $M\cap\{z=0\}$ consists
of the three components $X$,  $S$, and $\rho_X(S)$, each of which
is a smooth, embedded curve.
Now suppose that $\overline{S}$ has an endpoint $p$. Then  $p$
must  be a singular point of $M\cap\{z=0\}$. Thus the tangent
plane to $M$ at $p$ must be horizontal. By
statement~\eqref{claim:NoHorizontalPlanes} of  the claim, the
point $p$ must lie on $X\setminus\{0\}$. By the $\rho_Y$-symmetry
of $M$, the point $-p$ must also be an endpoint of $\overline{S}$.
Thus $S$  consists of $X$ together with the simple closed curve
$\overline{S}\cup \rho_X(S)$, which intersects $X$ precisely at
the the two points $p$ and $-p$.   This completes the proof of the
claim.
\end{proof}

We have established that $M$ is a complete, properly immersed
minimal surface with finite topology, one end, and bounded
curvature, such that the level set
\[
   M\cap \{z=0\}
\]
consists of finitely many curves with finitely many singular
points (points where curves cross). According to a theorem of
Rodriguez and Rosenberg~\cite{RodRos}, {\it any} minimal surface
$M$ with these properties has finite type, meaning:
\begin{enumerate}[   \upshape(a)]
 \item\label{FiniteType:ConformalType} $M$ is conformally a once-punctured riemann surface.
    Equivalently, the one-point compactification $M\cap\{\infty\}$ is conformally a
    compact riemann surface (in our case a torus).
 \item The one-form $(dg)/g$ (where $g$ is the Gauss map) is meromorphic  on
      $M\cup\{\infty\}$,
 \item Let $\eta$ be the holomorphic one-form on $M$ whose real part
 is $dv$, where $v(x,y,z)=z$ is the height function.  Then $\eta$ is a meromorphic one-form
 on $M\cup\{\infty\}$.
\end{enumerate}

In particular,~\eqref{FiniteType:ConformalType} is
assertion~(\ref{mainthm1:conformaltype}) of the theorem.

Since the height function $v$ is harmonic and nonconstant on $M$,
the one-form $\eta$ must have a pole at $\infty$ and nowhere else.
Since $\eta$ is a nonconstant meromorphic one-form on a torus, it
must have an equal number (counting multiplicity) of zeroes and
poles.  Thus it must have some zeroes on $M$.  In other words, the
height function $v$ must have critical points on $M$. By
statement~\eqref{claim:NoHorizontalPlanes} of the claim,  those
critical points must lie on $X\setminus\{0\}$.
 In particular, the level set $M\cap \{z=0\}$ is not everywhere smooth.
Thus by statement~\eqref{claim:Dichotomy} of the claim, the level
set consists of $X$ together with a simple closed curve that
intersects $X$ at two points.  This is
assertion~\eqref{mainthm1:level0}.

By elementary complex analysis, for each $c\in \RR$, the number of
ends of the level set $M\cap\{z=c\}$ (i.e., of $v^{-1}(c)$) is
equal to the order of the pole of $\eta$ at $\infty$.   For $c=0$,
there are two ends by assertion~\eqref{mainthm1:level0}. Thus
$\eta$ has a double pole at $\infty$, and the level set
$M\cap\{z=c\}$ has two ends at infinity for every $c$.

Now let $c\ne 0$.  Since the height function has no critical
points at height $c$, the level set
\begin{equation}\label{e:HalfofLevelc}
    M^+\cap \{z=c\}
\end{equation}
is a union of smooth embedded curves. None of the curves is closed
since $M^+$ is simply connected and embedded. Thus the level set
$M\cap \{z=c\}$ (which is obtained from~\eqref{e:HalfofLevelc} by
reflection in $Z$) is a disjoint union of non-closed smooth
curves.  We have just shown that this level set has exactly two
ends.   Thus it consists of exactly one curve.  This proves
assertion~\eqref{mainthm1:levelc}.

It remains only to show assertion~\eqref{mainthm1:asymptotic},
that $M$ is asymptotic to $H$ at infinity. This follows from
immediately from a theorem of Hauswirth, Perez, and Romon, who
prove that any embedded minimal surface of finite type, one end,
bounded curvature, and infinite total curvature must be asymptotic
to a helicoid at infinity ~\cite{HausPerRom}.

The Hauswirth-Perez-Romon Theorem is very general, but has a
rather lengthy proof. We can also deduce
assertion~\eqref{mainthm1:asymptotic} from the following theorem
(due to Hoffman and McCuan), which has a much shorter proof.

\begin{theorem*}\label{SymmetricEnd} \cite{HOMc03}
Let $E\subset \RR^3$ be a properly immersed, minimal annular end
that is conformally a punctured disk. Suppose that $ \frac{dg}{g}$
and $\eta$  both have double poles at the puncture and that $\eta$
has no residue at the puncture. If $E$ contains a vertical ray and
a horizontal ray, then a sub-end of $E$ is embedded and asymptotic
to a helicoid.
\end{theorem*}

We have already proved that the end of $M$ satisfies all the
hypotheses except for two: that $\eta$ has no residue at $\infty$,
and that $dg/g$ has a double pole at $\infty$. (We showed in 
proving assertion~\eqref{mainthm1:levelc} that $\eta$ has a double
pole at $\infty$.)

Since $\eta$ has no poles on $M$, it has no residue at $\infty$ by
Stokes' theorem.

Since $dg/g$ and $\eta$ are meromorphic one-forms on
$M\cup\{\infty\}$, their ratio is a meromorphic function and hence
has a limiting value at value at infinity:
\begin{equation}\label{e:LimitingValue}
   \lim_{p\to\infty} \frac{dg/g}{\eta}(p) = \xi \in \CC\cup\{\infty\}
\end{equation}
Recall (\cite{HoffmanKarcherSurvey}, p.~15) that the principal
curvatures at a point are $\pm \kappa$ where
\[
    \kappa = 4\left(|g| + \frac1{|g|} \right)^{-2} \left|  \frac{dg/g}{ \eta} \right|.
\]
\begin{comment}
(See \cite{os1}, 9.2.  In the notation there,
$dg=g'(\zeta)\,d\zeta$ and $\eta=fg\,d\zeta$, where $\zeta$ is a
local complex coordinate on $M$.)
\end{comment}
On $Z$, the tangent plane is vertical, so $|g|=1$ and
\begin{equation*}%\label{e:kappaonZ}
    \kappa =  \left| \dfrac{dg/g}{ \eta} \right|.
\end{equation*}
Thus by~\eqref{e:LimitingValue},
\begin{equation}\label{e:LimitingCurvature}
  \lim_{p\to\infty, \,p\in Z}\kappa(p) = |\xi|.
\end{equation}

Since $M$ has bounded principal curvatures, $\xi\ne \infty$.

Note that on any interval $I\subset \ZZ^+$ of length $4\pi$, the
tangent plane to $H$ turns through angle $4\pi$.  Since $M^+$ lies
on one side of $H$, this forces the tangent plane to $M$ to turn
through an angle at least $3\pi$ on $I$.  Thus there must be a
point in $I$ at which
\[
   \kappa  \ge \frac{3\pi}{4\pi} = \frac34.
\]
In particular, there is a sequence of such points in $Z$ tending
to $\infty$, so $\xi\ne 0$ by~\eqref{e:LimitingCurvature}.

Since $\xi$ is finite and nonzero, ~\eqref{e:LimitingValue}
implies that $dg/g$ and $\eta$ have poles of the same order at
$\infty$. Since $\eta$ has a double pole at  $\infty$, so does
$dg/g$.

All the hypotheses of the Hoffman-McCuan Theorem are satisfied, so
$M$ is asymptotic to $H$. This completes the proof of
Theorem~\ref{thm:main theorem 1}.
\end{proof}

\begin{theorem}\label{thm:main theorem 2}
For each $h>\pi/2$, there is a complete, properly embedded minimal
surface $N=N(h)$ in $\RR^3$ such that
\begin{enumerate}
 \item The intersection $N\cap H$ of $N$ with the helicoid $H$
 consists of the $z$-axis together with the
 horizontal lines
   $$
     H \cap \{z=nh\}, \,  n\in\ZZ.
   $$
 Furthermore, $N\setminus H$ consists of congruent, simply
 connected components.
 \item $N$ is invariant under the screw motion $\sigma_{2h}$.
 \item The portion $M$ of $N$ in the slab $|z|<h$ is bounded by the two horizontal
 lines $\rho_{\pm h}(X)$.
 \item The quotient $N/\sigma_{2h}$ has finite total curvature
 and is conformally a twice-punctured torus.
 \item $\Tan_pN$ converges to a horizontal plane as
 $\dist(p,Z)\to\infty$.
 \item The quotient $N/\sigma_{2h}$ is asymptotic to $H/\sigma_{2h}$ at infinity.
\end{enumerate}
\end{theorem}

\begin{proof} Fix an $h > \pi/2$ and
choose a sequence $R_i\to \infty$.  By Propositions~\ref{prop:aip}
and~\ref{existence3}, for sufficiently large $i$ there is a
$\rho_Y$-invariant minimal embedded disk $D_i$ in $H^+$ with
boundary $\Gamma(R_i,h)$ and with the annular intersection
 Let
$$
  M_i = \overline{D}_i\cup \rho_Z(\overline{D}_i)
$$
and let $N_i$ be the $\sigma_{2h}$-invariant surface obtained from
$M_i$ by repeated Schwarz reflections.  (Or, equivalently, let
$N_i=\cup_{n\in\ZZ}\sigma_{2hn}M_i$.)

By Theorem~\ref{thm:MainCurvatureEstimate}, the curvatures of the
$N_i$ are uniformly bounded, so (passing to a subsequence) we may
assume that the $N_i$ converge smoothly to a complete embedded
$\sigma_{2h}$-invariant minimal surface $N$.  It follows that the
$M_i$ converge smoothly to the closure of the surface
$$
   M=\{(x,y,z)\in N: |z|<h\}.
$$
Note that $N$ is orientable since it is properly embedded in
$\RR^3$.  If $L$ is a line in $N$, then the rotation $\rho_L:N\to
N$ is orientation-reversing.  Thus
$$
  \sigma_{2h}: N\to N
$$
is orientation-preserving since it is the product of two such
rotations (corresponding to $L=X$ and $L=H\cap\{z=h\}$). Therefore
$N/\sigma_{2h}$ is orientable.

Assertions (1), (2), and (3) follow immediately from the
construction.

Exactly as in the proof of Theorem~\ref{thm:main theorem 1}, $M$
is homeomorphic to $T\setminus \Delta$ where $T$ is a torus and
$\Delta$ is a closed disk in $T$.  Since $\overline{M}$ is $M$
together with its two boundary lines, $\overline{M}$ is
homeomorphic to the union of $T\setminus \Delta$ with two disjoint
arcs of $\partial \Delta$. If we identify those two arcs, either
the result is non-orientable or else it is topologically a twice
punctured torus. Since $N/\sigma_{2h}$ is orientable and since it
is the result of identifying the two boundary lines of
$\overline{M}$, it follows that $N/\sigma_{2h}$ is topologically a
twice-punctured torus.

Note that the total curvatures of the $M_i$ are uniformly bounded
by the Gauss-Bonnet theorem. It follows that $N/\sigma_{2h}$ has
finite total curvature, which implies by Huber's
theorem~\cite{hu1} (or by \cite{os1}, Theorem~9.1)) that it is
conformally a punctured Riemann surface. Thus it is a conformally
a twice-punctured torus.

Assertion~(5), i.e., that $\Tan_pN$ becomes horizontal as
$\dist(p,Z)\to\infty$, follows from Theorem~\ref{translates}.
(Note that case (3) of that theorem does not arise because
Theorem~\ref{thm:BoundOnGeodesic} gives a uniform bound on the
lengths of the closed geodesics in the $\overline{D}_i$.)

We now show that  the two horizontal ends are
helicoidal rather than planar. (See~\cite{mr93} for a discussion
of the possible end behavior of embedded periodic minimal
surfaces.) If they were planar, one end would correspond to the
plane $z=0$ and the other to the plane $z=h$ (since $N$ contains
horizontal lines at those heights).  But $\rho_Y$ preserves
orientation on $N$ and reverses orientation on the plane $z=0$.
Thus the ends cannot be planar.  
(Section~5 of \cite{hoffwhite2} gives another proof that the ends
of $N/\sigma_{2h}$ are  asymptotic to a helicoid.)

It remains only to show that $D$ is asymptotic to $H$ as $x^2+y^2\to \infty$
rather than to some other helicoid.

Recall that $\overline{D_i}$ has a unique embedded, closed geodesic.
Note that the geodesic divides $\overline{D_i}$ into two components, and
that the portions of $X^+$ and of $\sigma_h(X^+)$ in  $\overline{D}_i\setminus \{0\}$
belong to the same connected component.  (They are joined by a helical arc
in $\overline{D}_i$).   Thus by smooth convergence of $D_i\to D$ and of the
corresponding geodesics, we see that the simple closed geodesic
in $\overline{D}$ divides it into two components, one of which, $D^{upper}$,
contains $X^+$ and $\sigma_hX^+$.  Consequently $X^-$ and $\sigma_{-h}(X^-)$ belong
to the other component.

Let $R$ be large enough that the closed geodesic and all points of $\overline{D}$ 
with vertical tangent planes lie in the open cylinder
of radius $R$ about $Z$.
Then for each $r\ge R$, 
\[
   \overline{D^{upper}} \cap \{ x^2+y^2=r^2\}
\]
is a smooth compact curve, on each component of  which the angle function
\[
  \theta: H^+ \to \RR
\]
is monotonic.  
Thus  no component is a closed curve, and since there are exactly two endpoints,
namely $(r,0,0)$ and $\sigma_h(r,0,0)$, there is exactly one component.
Furthermore, that component can be parametrized by $\theta$ and therefore
written as:
\[
    ( r\cos\theta, r\sin \theta, f(r, \theta)) \qquad (\alpha \le \theta \le \beta)
\]
for suitable $\alpha$ and $\beta$.
Now $\theta=0$ on $X^+$ and $\theta= h$ on $\sigma_h(X^+)$, so
$\alpha=0$, $\beta=h$, $f(r,0)=0$, and $f(r,h)=h$.

Since this is true for all $r\ge R$, we have proved that outside of the cylinder of radius $R$
about $Z$, 
 $\overline{D^{upper}}$ may be written as
 \[
   \{ (r\cos\theta, r\sin\theta, f(r,\theta)): r\ge R, \quad 0\le \theta \le h \}
\]
where
\[
   f(r,0)=0, \qquad f(r,h) = h. \tag{*}
\]

Since $D$ is asymptotic to a helicoid, $\partial f/\partial \theta$ converges to a constant
as $r\to\infty$, namely the pitch of the helicoid.   Thus \thetag{*} implies that the constant
is $1$ and that
\[
   f(r,\theta) = \theta + o(r)
\]
which implies that $D^{upper}$ is asymptotic to $H$.
\end{proof}

\section{A compactness theorem}\label{section: compactness}

Let $\Nn$ be the class of all symmetric, embedded genus-one
helicoids, periodic and nonperiodic.  If $N\in \Nn$ is periodic,
let $h(N)$ be the smallest $h>0$ such that $N$ is
$\sigma_{2h}$-invariant. If $N$ is nonperiodic, let $h(N)=\infty$.

For each $h>\pi/2$, we have proved existence of an $N\in \Nn$ with
$h(N)=h$. Furthermore, there are no $N\in\Nn$ with $h(N)\le\pi/2$
(see~\cite{hoffwhite2}, Section 3). However, $\Nn$ may conceivably
contain examples that do not arise from our construction.
Nevertheless, the following two theorems apply to all symmetric,
embedded genus-one helicoids.

\begin{theorem}\label{thm:compactness}
If $\eta>\pi/2$, then the class
$$
 \Nn_\eta := \{ N\in \Nn: h(N)\ge \eta \}
$$
is compact with respect to smooth convergence on compact sets of
$\RR^3$.  Furthermore, each surface $N\in \Nn$ has all the
properties listed in Theorem~\ref{thm:main theorem 1} if
$h(N)=\infty$, or in Theorem~\ref{thm:main theorem 2} if
$h(N)<\infty$.
\end{theorem}

\begin{proof}
If $R=\infty$ and/or $h=\infty$, let us (by a slight abuse of
notation) interpret the expression ``$\partial D=\Gamma(R, h)$''
to mean ``$\partial D=\partial Q_{R, h}$ and $D$ is asymptotic to
$Q_{R,h}$ at infinity''. Note that if $R=\infty$, then $\partial
Q_{R,h}$ has no helical portions, but rather consists entirely of
horizontal rays together vertical segments or rays.  In another
paper (\cite{hoffwhite2}, Section~2), we will prove that if $N\in
\Nn$, then
$$
    D := \{ (x,y,z)\in N\cap H^+: |z|<h(N) \}
$$
is a disk with $\partial D=\Gamma(\infty\,, h(N))$.

We claim that {\em all} of the results in
Sections~\ref{section:Vertical} and~\ref{section:uniform_est} of
this paper remain true for any disk $D$ with $\partial
D=\Gamma(R,h)$, even when $R$ and/or $h$ is allowed to be
infinite.   By Remark~\ref{remark:deformed helical arcs},
Propositions~\ref{vertical thm 1} and~\ref{VtangentZ} are true for
$D\cap C$ for any sufficiently large solid cylinder $C$ about $Z$,
and thus those propositions are also true for $D$. No changes are
required in any of the other proofs. (In particular, one works
directly with $D$ rather than with the portion of $D$ in a solid
cylinder.) By Remark~\ref{remark:aip not needed},
Theorem~\ref{thm:BoundOnGeodesic} is true when $R_i\equiv \infty$
even without assuming the annular intersection property.

Consider a sequence $N_i\in\Nn_\eta$. By the discussion above, the
$N_i$ and the disks
$$
   D_i := \{ (x,y,z)\in N_i\cap H^+: |z|<h(N_i)\}
$$
satisfy all the estimates in Sections~\ref{section:Vertical}
and~\ref{section:uniform_est}.  In particular, the curvature
estimates imply that a subsequence of the $N_i$ converges to a
limit $N$, and Theorem~\ref{thm:BoundOnGeodesic} implies that the
sequence of points $D_i\cap Y$ is bounded. The proofs of
Theorem~\ref{thm:main theorem 1} and Theorem~\ref{thm:main theorem
2} then show that $N$ satisfies all the conclusions of
Theorem~\ref{thm:main theorem 1} (if $h(N)=\infty$) or
Theorem~\ref{thm:main theorem 2} (if $h(N)<\infty$), and so in
particular $N\in \Nn_\eta$.
\end{proof}

We end by pointing out that our bounds on points with vertical
tangent planes are uniform:

\begin{theorem}
For every $\eta>\pi/2$, there is an $R<\infty$ with the following
property.  If $N\in \Nn_\eta$ and if $p$ is a point in the
fundamental domain
$$
   M = \{ (x,y,z)\in N:  |z|< h(N)\}
$$
such that $\Tan_pN$ is vertical, then $p$ lies on the $z$-axis or
in the solid cylinder
$$
   \{(x,y,z): x^2+y^2\le R^2,  |z|<2\pi \}.
$$
\end{theorem}

\begin{proof} This follows immediately from
Propositions~\ref{translates} and~\ref{thm:BoundOnGeodesic}.
\end{proof}

\appendix
\section{Existence of unstable minimal disks}

Here we prove the minimax principle used in
Section~\ref{section:DiskConstruction}.

Let $U$ be an open subset of $\RR^3$.  We will call $U$ {\bf mean
convex} provided no smooth minimal surface in $\overline{U}$ whose
boundary is in $U$ can touch $\partial U$. We will call $U$ {\bf
strictly mean convex} provided no smooth minimal surface $M$ in
$\overline{U}$ can touch $\partial U$ except at its boundary
$\partial M$.  For example, in case that $\partial U$ is smooth,
$U$ is mean convex if and only if the mean curvature of $\partial
U$ (with respect to the inward-pointing normal) is everywhere
nonnegative, and $U$ is strictly mean convex if and only if the
mean curvature is positive on dense subset of $\partial U$.

Suppose $U$ has piecewise smooth boundary.  If each face has
nonnegative mean curvature with respect to the inward-pointing
normal, and if the interior angle along each edge is less than
$\pi$, then $U$ is mean convex.

\begin{theorem}\label{minimax theorem in appendix}
 Suppose $\Gamma$ is a smooth, simple closed curve in $\RR^3$ and
that $D_1$ and $D_2$ are two disjoint strictly stable smooth
embedded minimal disks in $\RR^3$ with $\partial D_1=\partial
D_2=\Gamma$. Suppose also that $D_1$ and $D_2$ meet transversely
along $\Gamma$.  Then $\Gamma$ bounds a weakly unstable, smoothly
embedded minimal disk in the region between $D_1$ and $D_2$.
\end{theorem}

\begin{proof}
Let $U$ be the region between $D_1$ and $D_2$. Enlarge $U$
slightly by pushing each $D_i$ outward by $\eps\, u_i$ times the
outward unit normal, where $u_i$ is the first eigenfunction of the
Jacobi operator on $D_i$, normalized to have maximum value $1$. By
strict stability, the resulting enlarged domain $U_\eps$ will be
strictly mean convex for all sufficiently small $\eps>0$. (See
Proposition~\ref{prop:expanding} below.)  Since the first
eigenfunction is positive everywhere, $U_\eps$ contains $U$
provided $\eps>0$.

Now we appeal to the following theorem:

\begin{theorem}\label{degree theory theorem}
  Suppose that $W$ is a mean convex region in $\RR^3$ and  that $V$
is a smooth open subset of $\partial W$. Then:
\begin{enumerate}
 \item A generic, smooth,
 simple closed curve $C$ in $V$ is
 non-critical in that $0$ is not an eigenvalue of the Jacobi
 operator of any smooth embedded disk in $W$ with boundary $C$.
 \item If $W$ is strictly mean convex and if $C$ is a non-critical
 smooth simple closed curve in $V$ that bounds a disk in $V$,
 then $C$ bounds finitely many smooth embedded minimal disks
 in $W$, and the number of even-index disks is exactly one more
 than the number of odd-index disks.
\end{enumerate}
\end{theorem}

This theorem is stated for strictly convex $W$ with smooth
boundary in Theorem 2.1 of \cite{wh4}. But the proof given there,
which is inspired by Tomi-Tromba \cite{tt1}, actually establishes
the more general result stated here.

To continue the proof of Theorem~\ref{minimax theorem in
appendix}, let $\Gamma'$ be such a generic curve in either of the
two smooth faces of $U_\eps$.  By strict stability of $D_i$
($i=1,2$) and by the implicit function theorem (cf. Theorem 3.1 of
\cite{White-Indiana87}), if $\Gamma'$ is sufficiently close to
$\Gamma$, then $\Gamma'$ will bound a smooth embedded minimal disk
$D_i'$ that is a slight perturbation of $D_i$. In particular,
$\Gamma'$ bounds at least two strictly stable disks in $U_\eps$.
Thus by Theorem~\ref{degree theory theorem}, it must bound at
least one minimal embedded disk $D'$ in $U_\eps$ that has odd
index and that is therefore strictly unstable.

Now take a subsequential limit of such $D'$ for a sequence of
$\Gamma'$s converging to $\Gamma$. The resulting surface is a
weakly unstable disk $D(\eps)$ in $U_\eps$. Now let $D$ be a
subsequential limit of the $D(\eps)$ as $\eps\to 0$.  Then $D$ is
a weakly unstable smooth embedded minimal disk in $\overline{U}$.
It is not equal to $D_1$ or to $D_2$ since they are both strictly
stable.  Hence by the strict maximum principle, $D$ cannot touch
$D_1$ or $D_2$.  That is, $D\subset U$.
\end{proof}

\begin{proposition}\label{prop:expanding}
  Let $D$ be a smooth, embedded, strictly stable, orientable minimal
surface in $\RR^3$. Let $\nu$ be a unit normal vectorfield on $D$.
Suppose $u$ is an eigenfunction corresponding to the first
eigenvalue of the Jacobi operator on $D$.  Then for all
sufficiently small $t\ne 0$, the surface $D(t)$ parameterized by
  $x\in D \mapsto x + t u(x)\nu(x)$
lies on one side of $D$ and has nonzero mean curvature vector that
points into the region between $D$ and $D(t)$.
\end{proposition}

\begin{proof}
Let $j$ be a positive integer and define an operator
 $\Phi: C^{j+2}(\overline{D}) \to C^j(\overline{D})$
on a neighborhood of the zero function as follows. Given $w\in
C^{j+2}(\overline{D})$, form the disk $x\mapsto x+w(x)\nu(x)$. Let
$\vec{h}(w)(x)$ be the mean curvature vector of this surface at
the point $x+w(x)\nu(x)$, and let
$$
   \Phi(w)(x) = \vec{h}(x)\cdot \nu(x).
$$
Then $\Phi$ is a smooth map of Banach spaces.  Also, $\Phi(0)=0$
and $\Phi'(0)$ is the Jacobi operator $J$.  Thus if $w$ is any
smooth function on $\overline{D}$, then
$$
   \Phi(w) = J(w) + O(|w|^2).
$$
In particular, letting $w=t u$ gives
$$
  \Phi(t u) = -t\lambda u + O(t^2)
$$
where $\lambda$ is the first eigenvalue of $J$.  By strict
stability, $\lambda>0$.  Thus
\begin{equation}\label{e:Cj convergence}
  \text{
   $\frac{\Phi(t u)}{t} \to -\lambda u$
   in $C^j$ as $t\to 0$. }
\end{equation}
Now $u$ does not vanish on the interior of $D$, and hence by the
boundary maximum principle $Du$ never vanishes at the boundary.
Thus the $C^j$ convergence \eqref{e:Cj convergence} implies (for
all sufficiently small $t\ne 0$) that $\Phi(t u)$ never vanishes
on the interior of $D$ and that it has the sign indicated in the
statement of the lemma.
\end{proof}

\begin{corollary}\label{cor:appendix}
Suppose $D$ is a strictly stable, embedded orientable minimal
surface in $\RR^3$. Then there is an open set $W$ containing $D$
with the following property.  If $M$ is a minimal surface in $W$
with $\partial M\subset D$, then $M\subset D$.
\end{corollary}

\begin{proof}  Choose $\delta>0$ so that for $|t|\le \delta$,
the disk $D(t)$ has the property asserted in
Proposition~\ref{prop:expanding}.
 Let
$$
    W = \bigcup_{|t|<\delta} D(t).
$$
Given a minimal surface $M$ in $W$ with $\partial M\subset D$,
 let $T$ be the largest number such that $M\cap D_T$ is
nonempty. Then $T=0$, since if $T$ were positive, the maximum
principle would be violated at the point of contact of $M$ and
$D_T$. Likewise the smallest $T$ such that $\overline{M}\cap D_T$
is nonempty is $0$.
\end{proof}

\begin{remark} Using the boundary maximum principle, one can show
that the corollary is also true for $M\subset W$ with $\partial
M\subset \overline{D}$.
\end{remark}

\bibliography{variation}

\begin{thebibliography}{HKW99}

\bibitem[CHM89]{chm2}
M.~Callahan, D.~Hoffman, and W.~H. Meeks, III.
\newblock Embedded minimal surfaces with an infinite number of ends.
\newblock {\em Inventiones Math.}, 96:459--505, 1989.

\bibitem[CM04a]{ColdingMinicozzi1}
Tobias~H. Colding and William~P. Minicozzi, II.
\newblock The space of embedded minimal surfaces of fixed genus in a
  3-manifold. {I}. {E}stimates off the axis for disks.
\newblock {\em Ann. of Math. (2)}, 160(1):27--68, 2004.

\bibitem[CM04b]{ColdingMinicozzi2}
Tobias~H. Colding and William~P. Minicozzi, II.
\newblock The space of embedded minimal surfaces of fixed genus in a
  3-manifold. {II}. {M}ulti-valued graphs in disks.
\newblock {\em Ann. of Math. (2)}, 160(1):69--92, 2004.

\bibitem[CM04c]{ColdingMinicozzi3}
Tobias~H. Colding and William~P. Minicozzi, II.
\newblock The space of embedded minimal surfaces of fixed genus in a
  3-manifold. {III}. {P}lanar domains.
\newblock {\em Ann. of Math. (2)}, 160(2):523--572, 2004.

\bibitem[CM04d]{ColdingMinicozzi4}
Tobias~H. Colding and William~P. Minicozzi, II.
\newblock The space of embedded minimal surfaces of fixed genus in a
  3-manifold. {IV}. {L}ocally simply connected.
\newblock {\em Ann. of Math. (2)}, 160(2):573--615, 2004.

\bibitem[FM05]{fm05}
Leonor Ferrer and Francisco Mart{\'{\i}}n.
\newblock Minimal surfaces with helicoidal ends.
\newblock {\em Math. Z.}, 250(4):807--839, 2005.

\bibitem[Fuj88]{fu1}
H.~Fujimoto.
\newblock On the number of exceptional values of the {G}auss maps of minimal
  surfaces.
\newblock {\em Journal of the Math. Society of Japan}, 40(2):235--247, 1988.

\bibitem[HK97]{HoffmanKarcherSurvey}
David Hoffman and Hermann Karcher.
\newblock Complete embedded minimal surfaces of finite total curvature.
\newblock In {\em Geometry, V}, volume~90 of {\em Encyclopaedia Math. Sci.},
  pages 5--93, 267--272. Springer, Berlin, 1997.

\bibitem[HKW93]{howe3}
D.~Hoffman, H.~Karcher, and F.~Wei.
\newblock The genus one helicoid and the minimal surfaces that led to its
  discovery.
\newblock In {\em Global Analysis and Modern Mathematics}. Publish or Perish
  Press, 1993.
\newblock K. Uhlenbeck, editor, p. 119--170.

\bibitem[HKW99]{hkw8}
D.~Hoffman, H.~Karcher, and F.~Wei.
\newblock The singly periodic genus-one helicoid.
\newblock {\em Comment. Math.Helv.}, 74:248--279, 1999.

\bibitem[HM88]{hm9}
D.~Hoffman and W.~H. Meeks, III.
\newblock A variational approach to the existence of complete embedded minimal
  surfaces.
\newblock {\em Duke Math. J.}, 57(3):877--893, 1988.

\bibitem[HM90]{HoffmanMeeksHalfspace}
D.~Hoffman and W.~H. Meeks, III.
\newblock The strong halfspace theorem for minimal surfaces.
\newblock {\em Invent. Math.}, 101(2):373--377, 1990.

\bibitem[HM03]{HOMc03}
D.~Hoffman and J.~McCuan.
\newblock Embedded minimal annular ends asymptotic to the helicoid.
\newblock {\em Comm. Analysis and Geometry}, 11(4), 2003.

\bibitem[HPR01]{HausPerRom}
Laurent Hauswirth, Joaqu{\'{\i}}n P{\'e}rez, and Pascal Romon.
\newblock Embedded minimal ends of finite type.
\newblock {\em Trans. Amer. Math. Soc.}, 353(4):1335--1370 (electronic), 2001.

\bibitem[Hub57]{hu1}
A.~Huber.
\newblock On subharmonic functions and differential geometry in the large.
\newblock {\em Commentarii Mathematici Helvetici}, 32:181--206, 1957.

\bibitem[HW]{hoffwhite2}
D.~Hoffman and B.~White.
\newblock The geometry of genus-one helicoids.
\newblock {\em Commentarii Mathematici Helvetici}.
\newblock To appear.

\bibitem[HW02]{howe5}
D.~Hoffman and F.~Wei.
\newblock Deforming the periodic genus-one helicoid.
\newblock {\em Experimental Mathematics}, 11(2):207--218, 2002.

\bibitem[MO90]{mos1}
X.~Mo and R.~Osserman.
\newblock On the {G}auss map and total curvature of complete minimal surfaces
  and an extension of {F}ujimoto's theorem.
\newblock {\em J. Differential Geom.}, 31(2):343--355, 1990.

\bibitem[MR93]{mr93}
William~H. Meeks, III and Harold Rosenberg.
\newblock The geometry of periodic minimal surfaces.
\newblock {\em Comment. Math. Helv.}, 68(4):538--578, 1993.

\bibitem[MR05]{MR01}
W.~H. Meeks, III and H.~Rosenberg.
\newblock The uniqueness of the helicoid.
\newblock {\em Annals of Math.}, 161:1--32, 2005.

\bibitem[MY82a]{my1}
W.~H. Meeks, III and S.~T. Yau.
\newblock The classical {P}lateau problem and the topology of three-dimensional
  manifolds.
\newblock {\em Topology}, 21(4):409--442, 1982.

\bibitem[MY82b]{my2}
W.~H. Meeks, III and S.~T. Yau.
\newblock The existence of embedded minimal surfaces and the problem of
  uniqueness.
\newblock {\em Math. Z.}, 179:151--168, 1982.

\bibitem[Nit89]{ni2}
J.~{C.}~{C.} Nitsche.
\newblock {\em Lectures on Minimal Surfaces}, volume~1.
\newblock Cambridge University Press, 1989.

\bibitem[Oss63]{oss}
Robert Osserman.
\newblock On complete minimal surfaces.
\newblock {\em Arch. Rational Mech. Anal.}, 13:392--404, 1963.

\bibitem[Oss86]{os1}
R.~Osserman.
\newblock {\em A Survey of Minimal Surfaces}.
\newblock Dover Publications, New York, 2nd edition, 1986.

\bibitem[Rad71]{RadoBook}
Tibor Rad{\'o}.
\newblock {\em On the problem of {P}lateau}.
\newblock Springer-Verlag, New York, 1971.
\newblock Subharmonic functions, Reprint.

\bibitem[RR98]{RodRos}
Lucio Rodr{\'{\i}}guez and Harold Rosenberg.
\newblock Minimal surfaces in {$\bold R\sp 3$} with one end and bounded
  curvature.
\newblock {\em Manuscripta Math.}, 96(1):3--7, 1998.

\bibitem[Sch66]{schneider1}
R.~Schneider.
\newblock A note on branch points of minimal surfaces.
\newblock {\em Proc. Amer. Math. Soc.}, 17:1254--7, 1966.

\bibitem[SER88]{er3}
R.~{S}a Earp and H.~Rosenberg.
\newblock On values of the {G}auss map of complete minimal surfaces in {${\bf
  R}^3$}.
\newblock {\em Comment. Math. Helvetici}, 63:579--586, 1988.

\bibitem[TT78]{tt1}
F.~Tomi and A.~J. Tromba.
\newblock Extreme curves bound an embedded minimal surface of disk type.
\newblock {\em Math. Z.}, 158:137--145, 1978.

\bibitem[Whi87a]{wh5}
B.~White.
\newblock Curvature estimates and compactness theorems in 3-manifolds for
  surfaces that are stationary for parametric elliptic functionals.
\newblock {\em Inventiones {M}ath.}, 88(2):243--256, 1987.

\bibitem[Whi87b]{White-Indiana87}
Brian White.
\newblock The space of {$m$}-dimensional surfaces that are stationary for a
  parametric elliptic functional.
\newblock {\em Indiana Univ. Math. J.}, 36(3):567--602, 1987.

\bibitem[Whi89]{wh4}
B.~White.
\newblock New applications of mapping degrees to minimal surface theory.
\newblock {\em Journal of Differential Geometry}, 29:143--162, 1989.

\bibitem[WHW]{wwh1}
M.~Weber, D.~Hoffman, and M.~Wolf.
\newblock An embedded genus-one helicoid.
\newblock {\em Annals of Mathematics}.
\newblock To appear.

\bibitem[Xav81]{xa1}
F.~Xavier.
\newblock The {G}auss map of a complete non-flat minimal surface cannot omit 7
  points of the sphere.
\newblock {\em Annals of Math.}, 113:211--214, 1981.

\end{thebibliography}
\bibliographystyle{alpha}

\end{document}